\documentclass[12pt,leqno]{article}
\usepackage{amssymb}
\usepackage{amsmath}

\def\la{\lambda}
\def\al{\alpha}
\def\om{\omega}
\def\Om{\Omega}
\def\ve{\varepsilon}

\def\ga{\gamma}

\def\span{{\rm span}\,}
\def\rank{{\rm rank}\,}
\def\maxrank{{\rm maxrank}\,}
\def\minrank{{\rm minrank}\,}
\def\range{{\rm range}\,}

\def\Re{{\rm Re}\,}
\def\Im{{\rm Im}\,}

\def\trans{\,{}^t\!}
\def\sgn{{\rm sgn}\,}
\def\dst{\displaystyle}
\def\Sym{{\rm Sym}\,}
\def\GL{{\rm GL}\,}
\def\tr{{\rm tr}\,}

\def\A{{\cal A}}
\def\B{{\cal B}}
\def\C{{\cal C}}
\def\D{{\cal D}}

\def\G{{\cal G}}

\def\S{{\cal S}}

\def\P{{\cal P}}
\def\N{{\cal N}}
\def\M{{\cal M}}

\def\R{{\cal R}}
\def\V{{\cal V}}
\def\W{{\cal W}}

\def\g{{\frak g}}

\def\de{\partial}

\def\xy{{(x,y)}}

\def\CC{{\mathbb  C}}

\def\HH{{\mathbb  H}}
\def\KK{{\mathbb  K}}
\def\NN{{\mathbb  N}}

\def\RR{{\mathbb  R}}

\def\11{{\mathbb  1}}

\def\1{{\bf 1}}

\def \trans{\,{}^t\!}

\def\be{\begin{enumerate}}
\def\ee{\end{enumerate}}
\def\noi{\noindent}
\def\qed{\smallskip\hfill Q.E.D.\medskip}

\newtheorem{theorem}{Theorem}[section]
\newtheorem{proposition}[theorem]{Proposition}

\newtheorem{cor}[theorem]{Corollary}
\newtheorem{lemma}[theorem]{Lemma}
\newtheorem{remark}[theorem]{Remark}
\newtheorem{remarks}[theorem]{Remarks}

\newtheorem{assumption}[theorem]{Standing Assumptions}

\begin{document}

\title{Local solvability of second order differential operators
with double characteristics I:\\
Necessary conditions}
\author{Detlef M\"uller\\[2mm]}
\date{}
\maketitle

\begin{abstract}
This is a the first in a series of  two articles devoted to the
question of local solvability of  doubly characteristic 
 second order differential operators. Consider doubly
characteristic  differential operators of the form 
\begin{equation} 
L=\sum^m_{j,k=1}\al_{jk}(x) X_j X_k+\,\mbox{lower order
terms}\, ,
\end{equation}
where the $X_j$ are smooth real vector fields and the $\al_{jk}$ are
smooth complex coefficients forming a symmetric matrix
$\A(x):=\{\al_{jk}(x)\}_{j,k}.$  We say that
$L$ is essentially dissipative at $x_0,$ if there is some
$\theta\in\RR$ such that
$e^{i\theta}L$ is dissipative at $x_0,$ in the sense that $\Re
\big(e^{i\theta}\A(x_0)\big)\ge 0.$  For a large class
of doubly characteristic operators
$L$ of this form, we show that a necessary condition for local
solvability at 
$x_0$ is essential dissipativity of $L$ at $x_0.$ By means of 
H\"ormander's  classical necessary condition for local
solvability, the proof is reduced to the following question:

Suppose that $Q_A$ and $Q_B$ are two real quadratic forms on a
finite dimensional  symplectic vector space, and let
$Q_C:=\{Q_A,Q_B\}$ be given by the Poisson bracket of $Q_A$ and
$Q_B.$ Then $Q_C$ is again a quadratic form, and we may ask: When
can we find a common zero of  $Q_A$ and $Q_B$ at which $Q_C$ does
not vanish?

The study of this question will occupy most of the paper, and the
answers may be of independent interest. 

In the second paper of this series, building on joint
 work with F.~Ricci, M.~ÊPeloso and others, we shall study local
solvability of  essential dissipative left-invariant operators of the above form on
Heisenberg groups, in a fairly comprehensive way. 
Various examples exhibiting a kind of exceptional behaviour from previous joint works,
e.g., with G.~Karadzhov, have shown that there is little hope for a complete
characterization of locally solvable operators on Heisenberg groups.  However, the
 scheme of what rules local solvability of second order operators on Heisenberg
groups in general   becomes evident from our work.
\footnote[1]{
{\it keywords:}  linear partial differential operator, local solvability, doubly
characteristic, real quadric,\\ Poisson bracket
 }

\end{abstract}

\section{Introduction}\label{introduction}

Consider a second-order differential operator of order $k$ with smooth coefficients 
$$L=\dst{\sum_{|\al |\le k}}c_{\alpha}(x)D^{\al}$$ 
on an open
subset $\Om$ of $\RR^n$, where
$D^{\alpha}:=\dst{\left(\frac{\de}{2\pi i\,\de
x_1}\right)}^{\al_1}\cdots \dst{\left(\frac{\de}{2\pi i\,\de x_n}
\right)^{\al_n}}$.
\smallskip

\noi $L$ is said to be {\it locally solvable} at $x_0\in\Om$, if
there exists an open neighborhood $U$ of $x_0$ such that the
equation $Lu=f$ admits a distributional solution $u\in\D'( U)$
for every $f\in C_0^{\infty}( U)$ (for a slightly more general
definition, see  \cite{hoermander3}).

Around 1956, Malgrange and Ehrenpreis proved that every constant coefficient operator is locally
solvable, and shortly later H. Lewy produced the following example of a nowhere solvable
operator on $\RR^3$:

$$ Z=X-iY\, , \ \text{where}\ \ X:=\frac{\de}{\de x}-
\frac y2\frac{\de}{\de u},\
Y:=\frac{\de}{\de y}+\frac x2\frac{\de}{\de u}\ .$$ 
Not quite incidentally, $Z$ is a left--invariant operator on a 
2--step nilpotent Lie group, the Heisenberg group $\HH_1.$ 

This example gave rise to an intensive study of so-called principal type operators, which
eventually led, most notably through the work of H\"ormander, Maslov, Egorov,
Nirenberg--Tr\`eves and Beals--Fefferman, to a complete solution of
the problem of local solvability of such operators (see \cite{hoermander3}). 

\medskip

Let us recall some notation.  Denote by $p_k(x,\xi):=
\sum_{|\al |=k}c_{\al}(x)\xi^{\al}$ 
the  principal symbol of $L.$ We shall consider $p_k$ as an
invariantly defined function on the reduced cotangent bundle  
$\C:=T^*\Om\setminus 0=\Om\times  (\RR^n\setminus\{0\})$ of $\Om.$

\noi Let us   denotes by $\pi_1$ the
base projection $\pi_1:T^*\Om\rightarrow\Om$, $(x,\xi)\mapsto x.$ $T^*\Om$ carries a 
canonical $1$-form, which, 
in the usual coordinates, is given by $\al=\sum_{j=1}^n \xi_j dx_j,$ so
 that $T^*\Om$
has a canonical symplectic structure, given by the $2$-form $d\al=\sum_{j=1}^n
d\xi_j\wedge dx_j.$ In particular, for any smooth real function
$a$ on $\Om$, its corresponding
Hamiltonian vector field $H_a$ is well-defined, and explicitly given 
by 
$$H_a:=\sum_{j=1}^n\left(
\frac{\de a}{\de \xi_j}\frac{\de}{\de x_j}-\frac{\de a}{\de x_j}\frac{\de}{\de \xi_j}
\right).$$
 If $\gamma$ is an integral curve of $H_a$, i.e., if $\frac d{dt}\gamma(t)=
H_a(\gamma(t)),$ then $a$ is constant along $\gamma$, and $\gamma$ is called a {\it null
bicharacteristic} of $a$, if $a$ vanishes along $\gamma$. Finally, the
(canonical) Poisson bracket of two smooth functions $a$ and $b$  on $T^*(\Om)$ is given
by 
$$\{a,b\}:=d\al(H_a,H_b)=H_a b=
\sum_{j=1}^n\left(
\frac{\de a}{\de \xi_j}\frac{\de b}{\de x_j}-\frac{\de a}{\de x_j}\frac{\de b}{\de
\xi_j}
\right).$$

Let 
$$\Sigma =\{p_k=0\}\subset\C$$
 denote the {\it characteristic variety} of $L$. $L$ is
said to be of {\it principal type}, if $D_\xi p_2$ does not vanish on $\Sigma$ (or,
more generally, if for every $\zeta\in\Sigma$ there is a real  number $\theta$ such
that
$d(\Re(e^{i\theta}p_k))(\zeta)$ and
$\al(\zeta)$ are non--proportional). 

In 1960, H\"ormander proved the following
fundamental result on non--existence of solutions (see \cite{hoermander-grundlehren}):

\begin{theorem}[H\"ormander]\label{Ho} Suppose there is some $\xi_0\in
\RR^n\setminus\{0\}$ such that 
$$a(x_0,\xi_0)=b(x_0,\xi_0)=0\quad \text{and} \ \{a,b\}(x_0,\xi_0)\ne 0,$$
 where
$a:=\Re p_k$ and $b:=\Im p_k.$ Then $L$ is not locally solvable at $x_0.$ 
\end{theorem}

A complete answer to the question of local solvability of principal type operators
$L$ was eventually given in terms of  the following condition
$(\P)$ of Nirenberg and Tr\`eves: 
\smallskip

$(\P).\ \ $ The function $\Im( e^{i\theta}p_k)$ does not take both positive and
negative values along a null--bicharacteristic $\gamma_\theta(t)$ of $\Re
(e^{i\theta}p_k),$ for any
$\theta \in\RR$. 
\smallskip

\noi In fact, $L$ of principal type is locally solvable at $x_0$ if and only if $(\P)$
holds over some neighborhood of $x_0$.
Notice that this is a condition solely on the principal symbol of $L$. 

\medskip

In this article, we shall consider second order differential operators $L$ with double
characteristics.
Let 
\[
\Sigma_2:=\{(x,\xi)\in\C :\quad D_{\xi}p_2 (x,\xi)=0\}\, .
\]
By Euler's identity, $\Sigma_2$ is contained in the characteristic
variety $\Sigma,$ and so
$\Sigma_2$ consists of the doubly characteristic cotangent vectors 
of $L$.\par

\medskip

If $p_2$ is real, and if we assume that $\Sigma_2$ is a submanifold
of codimension $m<n$ in $\C$, such that $\rank\, D^2_{\xi}p_2=m$
for every $\xi\in\Sigma_2$, then the following has been shown in
\cite{mueller-ricci-survey} :
\smallskip

Given any $x_0\in\pi_1
(\Sigma_2)$, there exist suitable linear coordinates near $x_0$
such that, in those new coordinates,
$p_2$ can be written in the form
\begin{equation}\label{9a}
p_2(x,\xi)= \trans (\xi'+E (x)\,\xi''){\A} (x)(\xi'+E
(x)\,\xi'')
\end{equation}
with respect to some splitting of coordinates $\xi=(\xi',\xi'')$,
$\xi'\in\RR^m$, $\xi''\in\RR^{n-m}.$ Here,   
${\A}(x)\in\Sym(m,\RR)$ is non-degenerate and $E(x)\in
M^{m\times(n-m)}(\RR)$, and both matrices vary smoothly in
$x$ (as usually,  $ \Sym (n,\KK)$ denotes the space of all symmetric $n\times
n$ matrices over the field $\KK,$ where $\KK$ will be either $\RR$ or $\CC.$)
\medskip

Motivated by this result, let us assume here that $L$ is a complex
coefficient differential operator with smooth coefficients, whose
principal symbol is given by \eqref{9a} where, however, now 
${\A}(x)\in\Sym(m,\CC)$ is a complex matrix. We then write
\[
{\A}(x)=A(x)+iB(x),\quad   x\in\Om\, ,
\]
with $A(x), B(x)\in\Sym(m,\RR).$
Notice that \eqref{9a} means that, up to a factor $(2\pi)^{-2},$
$L$ can be written as
\begin{equation}\label{9c}
L=\sum^m_{j,k=1}\al_{jk}(x) X_j X_k+\,\mbox{lower order
terms}\, ,
\end{equation}
where $X_j$ is the real vector field
\begin{equation}\label{9d}
X_j=\frac{\de}{\de x_j}+\sum^n_{l=m+1}E_{jl}(x)\frac{\de}{\de
x_l}\, ,\quad j=1,\ldots ,m\, ,
\end{equation}
with ${\A} (x)=\{\al_{jk}(x)\}_{j,k}$,
$E(x)=\{E_{jl}(x)\}_{j,l}$. 

It may be worth while mentioning   that every second order differential operator
$L,$ whose principal part is of the form 
\begin{equation}\label{9c'}
\sum^m_{j,k=1}\beta_{jk}(x) Y_j Y_k,
\end{equation} 
with smooth real vector fields $Y_j$ which are linearly independent at $x_0,$ can be
put into this form, locally near $x_0.$ This can easily be seen by means of a suitable
change of coordinates and choice of a suitable new basis $X_1,\ldots,X_m$ of the
$C^\infty$-module spanned by $Y_1,\ldots,Y_m.$ 

Moreover, by
\eqref{9a},
\begin{eqnarray*}
\Sigma_2 \supset  \{(x,\xi)\in\C : \xi'+E(x)\,\xi''=0\}\subset
          \{(x,\xi)\in\C : dp_2(x,\xi)=0\}\, ,
\end{eqnarray*}
so that in particular $\pi_1(\Sigma_2)=\Om,$ and equality holds in
these relations, if, e.g., $\A(x)$ is non-degenerate.  Finally,
denote by
$q_j$ the symbol of $X_j$ (up to a factor $2\pi$), i.e.,
\[
q_j(x,\xi):=\xi_j+\sum^n_{l=m+1}E_{jl}(x)\xi_{l}\, ,\quad
j=1,\ldots ,m\, ,
\]
and define the skew-symmetric matrix
\[
J_{(x,\xi'')}:=\big(\{q_j,q_k\}(x,\xi)\big)_{j,k=1,\ldots
,m}\;
\]
(compare \cite{mueller-ricci-survey}). Here, $\{\cdot,\cdot\}$ denotes again the
canonical  Poisson bracket on $T^*\Om$.  
Observe that
$J_{(x,\xi'')}$ depends indeed only on $x$ and the
$\xi''$-component of $\xi$. Notice also that if
$\A(x)$ is non-degenerate, then $J_{(x_0,\xi''_0)}$ is
non-degenerate if and only if
$\Sigma_2$ is symplectic in a neigborhood of $(x_0,\xi_0),$ where
$\xi_0'$ is chosen so that $(x_0,\xi_0)\in \Sigma_2$ (see, e.g. 
\cite{treves}, Proposition 3.1, Ch. VII).

If $J_{(x,\xi'')}$ is non-degenerate, then we can 
associate to $J_{(x,\xi'')}$ the
skew form
\[
\om_{(x,\xi'')}(v,w):= \trans v \,\trans (J_{(x,\xi'')})^{-1}\,
w,\quad v,w\in\RR^m,
\]
 which defines a symplectic structure on $\RR^m.$  In
particular, $m=2d$ is even.

Recall that if $\om$ is an arbitrary symplectic
form on $\RR^{2d},$ then we can associate to any smooth function $a$ on $\RR^{2d}$ the
Hamiltonian vector field $H^\om_a$  such that $\om(H^\om_a,Y)=da\cdot Y$ for all
vector fields $Y$ on
$\RR^{2d},$ and define the associated  Poisson bracket
accordingly by
\[
\{a,b\}_\om :=\om (H^\om_a,H^\om_b)\, .
\]

We can thus define a Poisson structure   $\{\cdot ,\cdot\}_{(x,\xi'')}$ on  $\RR^m$
(depending on the point $(x,\xi'')$) by putting $\{\cdot ,\cdot\}_{(x,\xi'')}:=\{\cdot
,\cdot\}_{\om_{(x,\xi'')}}.$

\medskip
In order to formulate our main theorem, we neeed to introduce some further notation
concerning quadratic forms.
\medskip

If ${A}\in\Sym(n,\KK),$ we shall denote by $Q_{A}$ the associated quadratic 
form 
$$Q_{A}(z):=\trans z{A} z,\quad  z\in\KK^n,$$ on $\KK^n.$ 
For any non-empty subset $M$  of a $\KK$- vector space
$V,$ 
 $\span_{\KK} M$ will denote its linear span over $\KK$ in $V.$

Let $A, B\in\Sym(m,\RR).$ We say that $A,B$  
form a {\it non-dissipative pair,} if $0$ is the only positive-semidefinite
element in $\span_\RR \{A,B\}$. Moreover, we put 
\begin{eqnarray*}
\maxrank \{A,B\} &:=&\max\{\rank F:F\in \span_\RR \{A,B\}\}\\
\minrank \{A,B\}&:=&\min\{\rank F:F\in \span_\RR \{A,B\}\,, F\ne 0\}.
\end{eqnarray*}

 Observe finally that if
$Q_A$ and
$Q_B$ are quadratic forms on $\RR^m$, then
$\{Q_A,Q_B\}_{(x,\xi'')}$ is again a quadratic form. 
\\

We can now state our main result.
\begin{theorem}\label{9E}
Let $L$ be given by \eqref{9c}, and let $x_0\in\Om$. Assume that
\begin{itemize}
\item[(a)]
$A(x_0), B(x_0)$ forms a non-dissipative pair.
\item[(b)]
There exists some
$\xi''_0\in\RR^{n-m}\setminus\{0\}$ such that $J_{(x_0,\xi''_0)}$
is non-degenerate, and the matrices  $A(x_0),B(x_0)$ and
$C(x_0,\xi''_0)$ are linearly independent over $\RR$, where
$C(x_0,\xi''_0)\in\Sym (m,\RR)$ is defined by 
\[
Q_{C(x_0,\xi''_0)}:=\{Q_{A(x_0)},Q_{B(x_0)}\}_{(x_0,\xi''_0)}\, .
\]

\item[(c)] Either
\be
\item[(i)]   $\minrank \{A(x_0), B(x_0)\}\ge 3$ and 
$\maxrank\{A(x_0), B(x_0)\}\ge 17, $ or
\item[(ii)] $\minrank \{A(x_0), B(x_0)\}= 2,$
$\maxrank\{A(x_0), B(x_0)\}\ge 9, $ and  the joint kernel 
 $\ker A(x_0)\cap \ker B(x_0)$ of $A(x_0)$ and $ B(x_0)$ is
either trivial, or a symplectic subspace with respect to the
symplectic form $\om_{(x_0,\xi_0'')}.$
\ee
\end{itemize}
 Then $L$ is not locally solvable at $x_0.$
\end{theorem}

Notice that, like condition  $(\P),$ the condition in (a) is again a sign
condition on the principal symbol of $L.$

\medskip
Theorem \ref{9E} shows that a "generic" operator $L$ of the form \eqref{9c} can be
locally solvable at $x_0$ only, if there is some $\theta\in\RR$ such that
$\Re(e^{i\theta}\A(x_0))\ge 0.$ 
A major remaining  task will thus be to study local solvability of $L$ under the
assumption that $A(x)=\Re \A(x)\ge 0$ for every $x\in \Om.$ A stronger condition is the
condition 
\begin{equation}\label{cone}|B(x)|\le A(x),\quad x\in\Om.
\end{equation}
This condition is equivalent to Sj\"ostrand's {\it cone condition} \cite{sjoestrand}.
It implies hypoellipticity with loss of one derivative of the transposed operator
$\trans L,$ for "generic" first  order terms in 
\eqref{9c}, and thus local solvability of $L$ at $x_0$ (see  \cite{hoermander3},
Ch. 22.4, for details and further references). 

Since, however, local solvability of $L$ is in general a much weaker condition than 
hypoellipticity  of $\trans L,$ we are still rather  far from understanding what rules
local solvability in general, even when the cone-condition is satisfied.

Nevertheless, for the case of homogeneous, left-invariant second order differential
operators on the Heisenberg group $\HH_n,$ a rather complete answer had been given in
\cite{mueller-ricci-cone}, and in the sequel \cite{mueller-pos} to the present
article, we shall extend these results by dropping the cone condition, thus  giving  a
fairly comprehensive answer  for left-invariant operators on Heisenberg groups. 

We should like to mention that, even if the cone-condition is satisfied, for instance  
small perturbations of the coefficients of the  first order  terms preserving the
values at
$x_0,$  may  influence local solvability and lead to local solvability in
situations where the unperturbed operator is not locally solvable at $x_0$ (see, e.g.,
\cite {christ-karadzhov-mueller}). Moreover,  if, e.g., $\maxrank\{A(x_0), B(x_0)\}=5$
or
$7,$ then the statement in Theorem \ref{9E} are known to be wrong (see
\cite{mueller-karadzhov},\cite{mueller-peloso-nonsolv}).

All these results indicate that there is rather little hope for a complete
characterization of local solvability for doubly characteristic operators in general,
but that Theorem \ref{9E} in combination with the above mentioned results on
hypoellipticity give at least rather satisfactory  answers in the  "generic" case.

\medskip

 Theorem \ref{9E} can be reduced by means of H\"ormander's Theorem \ref{Ho} to the
following  result concerning real  quadrics, which may also be of independent 
interest and which represents the core of this work. 

\begin{theorem}\label{8ZZ}
Assume that $\RR^n=\RR^{2d}$ is endowed
with the canonical symplectic form. Let $A,B\in \Sym(n,\RR)$ be
linearly independent, and assume that $A,B$ forms a
non-dissipative pair.  Define $Q_C:=\{Q_A,Q_B\}$ as the Poisson bracket of $Q_A$
and $Q_B,$ and assume that $A,B$ and $C$ are linearly independent.

Then there exists a point $x\in\RR^n$
such that 
$$Q_A(x)=Q_B(x)=0 \ \text{and}\ Q_C(x)\ne 0,$$ 
provided one of the following conditions are satisfied:
\be
\item[(i)]   $\minrank \{A,B\}\ge 3$ and 
$\maxrank\{A,B\}\ge 17, $
\item[(ii)] $\minrank \{A,B\}= 2,$
$\maxrank\{A,B\}\ge 9, $ and the joint radical
$\R_{A,B}:= \ker A\cap\ker B$ of $Q_A$ and $Q_B$ is either trivial,
i.e., $\ker A\cap \ker B=\{0\},$  or a symplectic subspace of
$\RR^n.$ 
\ee

\end{theorem}
\bigskip

The article is organized as follows. Sections 2 and 3 are devoted to the proof of
Theorem \ref{8ZZ}. Notice that  this theorem essentially states  that
the quadratic form
$Q_C:=\{Q_A,Q_B\}$ can only vanish on the joint zero set $\{Q_A=0\}\cap \{Q_B=0\}$ of
two  linearly independent quadratic forms $Q_A$ and $Q_B$  forming   a
non-dissipative pair, if $C$ is a linear combination of $A$ and $B.$ 

Of course, this
can only be true if $\{Q_A=0\}\cap \{Q_B=0\}$ is sufficiently big, and we shall show
in Section 2 that  (under these assumptions )
the quadrics  $\{Q_A=0\}$ and $ \{Q_B=0\}$ do in fact intersect tranversally in a
variety $\N$ of dimension $n-2.$ 

In Section 3.1, we prove some auxiliary results and recall some basic notions and
facts on semi-algebraic sets. 

The proof of Theorem \ref{8ZZ} is then given in Sections 3.1 and 3.2. We distinguish
between the situation where no stratum of $\N$ spans $\RR^n$ (Section 3.1) and the
case where at least one stratum spans (Section 3.2). It is interesting to notice that
the condition that $Q_C$ be  the Poisson bracket of $Q_A$ and $Q_B$ is only needed in
the first case (see Theorems \ref{8Z} and \ref{MF2}) . We also present a number of
examples in order to demonstrate that the conditions in the main Theorem \ref{8Z} of
Section 3.1 are essentially necessary. 

Section 4 finally contains the proof of Theorem \ref{9E}. Moreover, we give various
applications of this theorem to left-invariant differential operators on
2-step nilpotent Lie groups (compare Corollary \ref{9j} for general $2$-step nilpotent
Lie groups, and Corollary \ref{9k} for the particular case of the Heisenberg group). We
also indicate that Theorem \ref{9E} has applications to higher step situations too,
for instance on  $r$-step nilpotent Lie groups.

\setcounter{equation}{0}
\section{On the intersection of two real quadrics}\label{formproblem}

If $V$ is a  $\KK$- vector space,
and if $v_1,\dots, v_k$ are vectors in $V,$ then
$v_1\wedge\cdots\wedge v_k$ will denote their exterior product
in $\Lambda^k(V).$ In particular, $v_1,\dots, v_k$ are linearly
 dependent if and only if $v_1\wedge\cdots\wedge v_k=0.$ The
open interior of a subset S of some topological space will be denoted 
by $S^0.$

If $M$ is a non-empty subset of $\Sym(n,\RR),$ then we say that $M$
is {\it non-dissipative}, if $0$ is the only positive-semidefinite
element in $\span_\RR M$.          

\begin{lemma}\label{8A} {Let $M$ be a non-empty subset of
$\Sym (n,\RR)$}. Then the following are equivalent: %%
\be %%
\item[(i)] $M$ is non-dissipative.%%
\item[(ii)] There is some positive definite matrix $ Q>0$
such that 
\begin{equation}
            \tr (\trans Q F Q)=0\quad \text{for every}\ F\in M.
\label{tr0}
\end{equation}
\item[(iii)] There is a matrix $T\in \GL(n,\RR)$ 
such that 
\begin{equation}
            \tr (\trans T F T)=0\quad \text{for every}\ F\in M.
\label{tr1}
\end{equation}
\ee
\end{lemma}

\noi{\bf Proof.} (i) $\Rightarrow$ (ii). Let $V=\Sym(n,\RR)$, and
let $\mathcal{P}\subset V$ denote  the closed cone of positive
semidefinite  matrices in $V.$ Put $K:=\{E\in\P: \tr E=1\}.$ %%
Notice that  $E\in \P$ has vanishing trace if and only if $E=0.$
Therefore, $\P\setminus \{0\}=\bigcup_{t>0} tK.$ 

If $W:=\span_\RR M,$ then $K$ and $W$ are  convex 
subsets of $V,$ which are disjoint, by (i). Moreover, $K$ is compact and 
$W$ is closed. By Hahn-Banach's theorem (see, e.g. \cite{rudin},
Theorem 3.4 (b)), there exists a linear functional 
$\mu\in V^*$ and $\ga\in\RR$, such that 
\[
\mu(E)>\ga >\mu (F)\qquad  \text{for all}\ \; E\in K,\, F\in W\, .
\]
Since $W$ is a linear space, this implies $\mu|_W =0, $ hence $\mu(E)\ge
0 \  \text{for all}\  E\in K.$ 

Choose $P\in V$ such that $\mu(M)=\tr(PM)\ \text{for all}\ \; M\in V.$
Then 
$$\tr(PE)>0\quad \text{for all}\ \;E\in\P.$$
Rotating coordinates, if necessary, we may assume that $P$ is 
diagonal, say $P=\mathrm{diag}(\la_j).$ But then clearly 
$\la_j>0,\, j=1,\ldots ,n$, hence $P>0$. Choose $Q>0$ such that 
$P=Q^2$. Then%%
\[%%
0=\tr(Q^2F)=\tr(^tQFQ)\quad  \text{for all}\ \; F\in W\, ,
\]%%
so that \eqref{tr0}Ê holds. 
\smallskip

\noi (ii) $\Rightarrow$ (iii) is trivial.

\smallskip

\noi (iii) $\Rightarrow$ (i). Assume that $F\in W$ and $F\ge 0.$
Then (iii) implies that $\tr(\trans TFT)=0$, where $\trans TFT\ge 0$.
 This implies $\trans TFT=0$, hence $F=0$.

\qed

If $A\in\Sym(n,\RR),$ then we put 
\[
\Gamma_{{A}}:=\{z\in\RR^n:Q_{{A}}(z)\le 0\}\, .
\]

Recall that a pair ${A},{B}\in\Sym (n,\RR)$ such that $\{{A},{B}\}$ is
non-dissipative  is called a  non-dissipative pair.  Notice
that this property depends only on the linear span of
${A}$ and
${B}.$ In view of Lemma \ref{8A}, it will sometimes be convenient to
assume that a linear change of coordinates has been performed so
that ${A},{B}$ have vanishing trace. 

%\medskip
 
%In the sequel, we may therefore assume w.r. that
%\begin{equation}\label{8a}
%tr A=tr B=0
%\end{equation}

%\noi (replacing $A$ by $\trans QAQ$ etc.) Moreover, we may assume
%$A=\mathrm{diag}(\la_j)$, with $\dst{\sum_j}\la_j=0$.\\

%\begin{remark}\label{8B} Since we only use $C$ independently of
%$A$ and $B$, the above reduction is justified (we do not need that
%$Q\in \mathrm{Sp}(u,\RR!)$.
%\end{remark}

By $S^{n-1}$ we shall denote the
Euclidean unit sphere in
$\RR^n,$ and by $B_r(x)$ the open Euclidean ball of radius $r$ centered
at $x\in\RR^n.$

\begin{theorem}\label{8C}
Let ${A},{B}\in\Sym(n,\RR),$ and assume that
\begin{equation}\label{8b}
\tr{A}=\tr{B}=0
\end{equation}
and
\begin{equation}\label{8c}
\Gamma_{{A}}\subset\Gamma_{{B}}.
\end{equation}
Then there is some $c\in\RR$ such that 
\begin{equation}\label{8d}
{B}=c{A}\, .
\end{equation}
\end{theorem}

\noi {\bf Proof.} After a rotation of coordinates, we may assume that
\[
{A} =\left(\begin{array}{cc}A_1 & 0\\
                        0   & -A_2  \end{array}\right)\, ,
\]
w.r. to the decomposition $\RR^n=\RR^k\times\RR^\ell$, with
(diagonal) matrices $A_1>0,A_2\ge 0$. Write correspondingly
\[
{B}=\left(\begin{array}{cc} B_1   & B_3\\
                           ^tB_3 & -B_2\end{array}\right)\, .
\]
%with $B_1,B_2>0$ (if ${B}$ is sufficiently close to ${A}$, this can
%be achieved).
 We decompose $z\in\RR^n$ as
$z=(x,y)\in\RR^k\times\RR^\ell$. 

The case ${B}=0$ is trivial, so let us assume that ${B}\ne 0.$ Observe
first that 
$$B_2\ge 0,$$
for, if $y\in\RR^\ell,$ then  $(0,y)\in\Gamma_{A},$  hence $-\trans
yB_2 y\le 0,$ by \eqref{8c}.

We claim that, for every $y\in\RR^\ell,$ 
\begin{equation}\label{8pos}
\al(y):=\trans yA_2 y>0 \quad \implies \quadÊ\beta(y):=\trans yB_2 y>0.
\end{equation}
Indeed, if $\Om:=\{y\in\RR^\ell: \al(y)>0\},$ and if $y\in\Om,$ then,
given $x\in\RR^k,$ there is some $r_x>0$ such that
$(x,ty)\in\Gamma_{A}$ whenever $t\in\RR,|t|>r_x.$ We thus find that 
$$\trans x B_1x+ 2t(B_3y)\cdot x -t^2 \beta(y)\le 0 \quad
\text{provided}\  |t|>r_x.$$
If $\beta(y)=0,$ choosing both signs of $t,$ we see that this implies 
$$\trans x B_1x\le 0 \quad\text{and}\  (B_3y)\cdot x=0 \quad\text{for
all}Ê\ x\in\RR^k,$$
hence $-B_1\ge 0$ and $B_3y=0.$ Therefore, the matrix 
\[
-\left(\begin{array}{cc} B_1   & 0\\
                           0 & -B_2\end{array}\right)
\]
is postive semi-definite. Since it has vanishing trace, it must
vanish, so that $B_1=0,B_2=0,$ and $B_3y=0.$ 

Thus, if $\beta(y_0)=0$ for some $y_0\in\Om,$ then $\beta(y)=0$  and
consequently  $B_3y=0,$  for every $y\in\Om.$ Since $\Om$ is dense in
$\RR^\ell,$ we obtain $B_3=0, $ hence ${B}=0,$  contradicting our
assumption on ${B}.$ This proves \eqref{8pos}.

\noi For $e\in S_\Om:=S^{\ell-1}\cap \Om\subset\RR^\ell$, put
\[
Q^e_{{A}}(x):=Q_{{A}}(x,e)=\trans x A_1x-\al (e)\, ,
\]
where $\al (e)= \trans eA_2 e>0\, .$
Similarly, let
\[
Q^e_{{B}}(x):=Q_{{B}}(x,e)= \trans xB_1x+2 (B_3 e)\cdot x-\beta (e)\, ,
\]
where $\beta (e)=  \trans e B_2 e>0 ,$ because of \eqref{8pos}.

Finally, for $e\in S_\Om$ fixed, let
\begin{eqnarray*}
f(x) & := & \frac 1 {\al (e)} Q^e_{{A}}(x)= \trans x\frac{A_1}{\al
(e)}x-1\, ,\\
g(x) & := & \frac 1{\beta (e)}Q^e_{{B}}(x)= \trans x\frac{B_1}{\beta
(e)}x +\xi\cdot x-1\, ,
\end{eqnarray*}
where $\xi:=2\frac{B_3e}{\beta (e)}$. From \eqref{8c}, we know that
\begin{equation}\label{8e}
f(x)\le 0\quad\implies\quad g(x)\le 0\quad \text{for all}\  x\in\RR^k.
\end{equation}

\begin{lemma}\label{8F}
Let $\A,\B\in\Sym(k,\RR),$ and assume that $\A>0.$  Moreover, let
$\xi\in\RR^k,$ and put%%
\begin{eqnarray*}
f(x) & := & \trans x \A x-1,\\
g(x) & := & \trans x \B x+\xi\cdot x-1\, .
\end{eqnarray*}
Then \eqref{8e} implies%%
\begin{equation}\label{8g}
\tr \A\ge \tr \B.
\end{equation}
Moreover, either $f=g$, or%%
\begin{equation}\label{8h}
\tr \A > \tr \B\, .
\end{equation}
\end{lemma}

\noi \noi {\bf Proof.} Observe that%%
\begin{equation}\label{8i}
\Delta f=2\, \tr \A,\quad \Delta g= 2\, \tr \B\, .
\end{equation}

 Assume now that $\tr \A\le \tr \B$.\par%%
\noi Then, by \eqref{8i}, $\Delta (g-f)\ge 0$, so that $g-f$ is
subharmonic. Moreover, $g(x)-f(x)\le 0$ for $f(x)=0$. By the maximum
principle, we thus conclude that
\[
g(x)-f(x)\le 0\;\mbox{in the ellipsoid}\; \{f\le 0\}\, .
\]

\noi Thus, there is some $\ve >0$, so that
\[
\trans x(\B-\A)x+\xi\cdot x\;\le\; 0\quad \text{for all}\ x\in
B_{\ve}(0)\, .
\]

\noi Then also $\trans x(\B-\A)x-\xi\cdot x\le 0,$ hence
\[
\trans x(\B-\A)x\; \le\; 0\quad \text{for all}\ x\in B_{\ve}(0)\, .
\]

\noi By homogeneity, this implies $\B-\A\le 0$, hence
\begin{equation}\label{8j}
\B\le \A\, .
\end{equation}

\noi In particular, $\tr \B\le \tr \A$, hence $\tr \A= \tr \B$.
But then $\A-\B\ge 0$, $\tr (\A-\B)=0,$ hence  $\A-\B=0$, so that 
$\A=\B$.\par%%
So, either $\tr \A> \tr \B$, or $\A=\B$. This proves \eqref{8g}.\par%%
\noi Moreover if $\A=\B$, then $g(x)=f(x)+\xi\cdot x$, and, after a
linear change of coordinates, we may assume that $\A=I$, i.e.,
\[
f(x)=|x|^2-1\, .
\]

\noi Thus, $|x|\le 1$ implies $|x|^2+\xi\cdot x\le 1$. If $\xi\not=
0$, choosing $x=\frac{\xi}{|\xi |}$, we obtain $1+|\xi |\le 1$, hence
$\xi=0$, a contradiction. Therefore, $\xi=0$, hence $f=g$.\\

\smallskip
 
\hfill Q.E.D.
\medskip

\noi Going back to the proof of Theorem \ref{8C}, we can now conclude
that
\[
\frac{1}{\al (e)} \tr A_1\;\ge\;\frac{1}{\beta (e)} \tr B_1\, ,
\]
\noi i.e., %%
\begin{equation}\label{8k}
\beta (e) \tr A_1\;\ge\;\al (e) \tr B_1\quad \text{for all}\ e\in
S_\Om\, ,
\end{equation}
and this  inequality is strict, unless $f=g$, i.e.,%%
\begin{equation}\label{8l}
\frac{A_1}{\al (e)}=\frac{B_1}{\beta (e)}\ \mathrm{and}\
\xi=(2/\beta(e))B_3 e=0\, .
\end{equation}
Notice that, by continuity and since $\S_\Om$ is dense in
$S^{\ell-1},$ \eqref{8k} holds indeed for all $e\in S^{\ell-1}.$

We distinguish therefore two cases. %\\
\medskip

\noi (a) If there exists some $e\in S_\Om$ such that \eqref{8k}
holds strictly, we choose an orthonormal basis $e_1=e,e_2,\ldots
,e_\ell$ of $\RR^\ell$. Then %%
\[
\beta (e_j)\,\tr A_1 \ge \al (e_j)\,\tr B_1,\quad j=1,\ldots, \ell\, ,
\]%%
and the inequality is strict for $j=1$. Summing in $j$, we thus
obtain%%
\[
\tr B_2\,\tr A_1 > \tr A_2\,\tr B_1\, .
\]\
But, since $0=\tr{A} = \tr{B}$, we have $\tr A_2=\tr A_1$, $\tr B_2=\tr
B_1$, hence %%
\[ %%
\tr B_1\, \tr A_1 > \tr A_1 \,\tr B_1\, ,
\]
a contradiction.%%
\medskip

\noi (b) There  remains the case where %%
\[%%
\beta(e)A_1=\al(e)B_1\ \text { and } B_3 e=0\quad
\text{for all}\  e\in S_\Om\, .%%
\]%%
Again, by continuity, this then holds for all $e\in S^{\ell-1}.$ 
But then $B_3=0$, hence%%
\[%%
{B}=\left( \begin{array}{cc} B_1 & 0\\
                            0   & -B_2\end{array}\right)\, .%%
\]
Moreover, $B_1=c A_1$ for some $c>0$. Then we see that%%
\[%%
\beta(e)=c\al(e)\, ,\quad \text{for all}\  e\in S^{\ell-1}\, ,%%
\]
which implies, by homogeneity,%%
\[%%
 \trans y B_2y= \trans y (cA_2)y\quad \text{for all}\  y\in\RR^{\ell}\,
,%%
\]
i.e., $B_2=cA_2.$ We thus get%%
\[%%
{B}=c{A}\, .%%
\]
\qed

\begin{proposition}\label{8M}
Let $E$ and $D$ be open ellipsoids in $\RR^n, n\ge 2$, whose
boundaries don't intersect transversally anywhere. Then, if $E\cap
D\not=\emptyset$, either  $E\subset D$ or $D\subset E$.
\end{proposition}

\noi {\bf Proof.} {\sl 1. Case $n=2$}.\\%%
We may then assume that the boundary  $\de D$ of $D$ is a circle
\[
\de D=\{(x,y)\in\RR^2:(x-\xi)^2+(y-\eta)^2=r^2\}\, ,
\]
and that
\[
\de E=\left\{(x,y)\in\RR^2:\frac{x^2}{a^2}+\frac{y^2}{b^2}=1\right\}\, ,
\]
where $a\ge b>0$.

Assume also that neither $E\subset D$ nor $D\subset E$.
Then, by convexity, $\de E$ can neither be contained in $\overline
D$, nor in $(\RR^2\setminus D)$, i.e., there is a point in $\de E$
lying in $\RR^2\setminus D$, and another one in $D$, hence, by
continuity of the boundary curve of $D$, $\de E\cap\de D\not=
\emptyset$.\par%%

Say that $E$ and $D$ {\sl pierce} at $X\in\de E\cap\de D$, if every 
neighborhood $U$ of $X$ contains a point in $\de E\setminus D$
and one in $\de D\setminus E$.\par

If $E$ and $D$ don't pierce at any $X\in\de E\cap\de D$, again by
continuity of the boundary curves of $E$ and $D$, either $\de
E\subset {\overline D}$ or $\de D\subset {\overline E}$, hence 
$E\subset D$ or $D\subset E$, by convexity.\par%%

Consequently, $E$ and $D$ pierce at at least one point $X_0\in\de
E\cap\de D$. By symmetry, we may assume that $X_0=(x_0,y_0)$, with
$x_0\ge 0$ and $y_0\ge 0$.\par

Assume first $0<x_0<a$, hence $0<y_0<b$.
If $D$ and $E$ pierce at $X$, then $\de D$ and $\de E$ must
have the same curvature at $X,$ namely $1/r,$ as can easily bee seen from 
the local Taylor expansions of the boundaries curves. If $a=b$, then this
implies
$D=E$, so assume $a>b$. Since $\de D$ has constant curvature, there are
exactly four points on $\de E$ which are potential piercing
points, namely $X_0,X_1:=(x_0,-y_0), X_2:=(-x_0,y_0)$ and
$X_4:=(-x_0,-y_0)$.
Moreover, by continuity, there must be at least one more piercing point,
besides $X_0$.\par

Assume, e.g., that $X_1$ is a second piercing point. Since
$\nu_0:=\left(\frac{x_0}{a^2},\frac{y_0}{b^2}\right)$ is normal to
$\de E$ at $X_0$, and
$\nu_1:=\left(\frac{x_0}{a^2},-\frac{y_0}{b^2} \right)$ is normal
to $\de E$ at $X_1$, and since the lines $X_0+\RR\nu_0$ and
$X_1+\RR\nu_1$ meet exactly at $(\xi
,\eta):=\left(\left(1-\frac{b^2}{a^2} \right)x_0,0\right)$, we see
that $\de D$ must be the circle with center $(\xi,\eta)$ and
radius $r:=\left(\left(\frac b a
\right)^4x^2_0+y^2_0\right)^{1/2}$.
Since $y^2_0=b^2-\frac{b^2}{a^2}x^2_0$, we thus have 
\[
r^2=b^2+\left( \frac{b^4}{a^4}-\frac{b^2}{a^2}\right)x^2_0\, .
\]

Moreover, computing the curvature of the ellipse $\de E$ at $X_0,$ and
comparing it with that of the circle $\de D,$  we then  find that 
\[
\frac{b/a^2}{\left (1+(\frac {b^2}{a^2}-1)\frac{x_0^2}{a^2}\right )^{3/2}}
=\frac{1/b}{\left (1+(\frac {b^2}{a^2}-1)\frac{x_0^2}{a^2}\right
)^{1/2}},
\]
hence $
0=(1-\frac {b^2}{a^2})(1-\frac {x_0^2}{a^2}).
$
This implies $x_0=a,$ in contradiction to our assumptions.
\medskip

The case where $X_3$ is a second piercing point can be treated in
a similar way. And, $X_4$ cannot be a piercing point, since
$X_4=-X_0$, so that the center of $\de D$ would have to ly on the
two parallel lines $X_0+\RR\nu_0$ and $-X_0+\RR (-\nu_0)$, which
would imply $\nu_0=sX_0$ for some $s\in\RR$. But this is impossible, since
$a>b$.\par%%

The cases where $x_0=0$ or $x_0=a$ are even easier, as well as the
case where $a=b$, and are left to the reader.
\bigskip

\noi {\sl 2. Case $n\ge 3$.}\\
\noi If neither $E\subset D$ nor $D\subset E$, we may choose points
$x_0\in E\cap D$, $x_1\in E\setminus D$ and $x_2\in D\setminus 
E$ spanning a two-dimensional affine plane $V.$  We may then reduce the
problem to the 2-dimensional case by restricting ourselves to the affine
plane $V.$ \par%%

Indeed, if $E:=\{P<0\}$, $D:=\{Q<0\}$, for suitable elliptic
quadratic functions $P$ and $Q$, then let
\begin{eqnarray*}
p(s, t) & := & P(x_0+s(x_1-x_0)+t(x_2-x_0)) = P\circ \gamma
(s, t)\, ,\\
q(s, t) & := & Q(x_0+s(x_1-x_0)+t(x_2-x_0)) = Q\circ \gamma
(s,t),\quad (s,t)\in\RR^2.
\end{eqnarray*}
Then $p$ and $q$ are elliptic quadratic functions on $\RR^2$, and
\[
E\cap V = \{p<0\}=: \tilde{E}, \  D\cap V =
\{q<0\}=:\tilde{D}\, ,
\]
in the coordinates $(s,t)$ for $V$.\par%%

And, if $p(s,t)=0=q(s, t)$, then $\nabla
p(s,t)\wedge\nabla q(s,t)=0$. For, if
$\nabla p(s,t)$ and $\nabla q(s,t)$
were linearly independent, then for $X=\gamma (s,t)\in\de
E\cap\de D,$ the vectors 
\begin{eqnarray*}
\nabla P(X)\cdot D\gamma (s,t) & = & 
\nabla p(s,t)\\
\nabla Q(X)\cdot D\gamma(s,t) & = & \nabla q(s,t),
\end{eqnarray*}
would be  linearly independent, hence also $\nabla
P(X)\wedge\nabla Q (X)\not= 0$.\par%%
  
Thus, we could apply Case 1 to $\tilde{E}$ and $\tilde{D},$ and 
conclude that one is contained in the other, say, e.g.,
$\tilde{E}\subset\tilde{D}$. But this would contradict our
assumption that $x_1\in E\cap V\setminus D\cap V$.\par%%

\qed

We shall need the following modification of Proposition \ref{8M}.

\begin{proposition}\label{8M0}
Let $A,B\in\Sym(n,\RR),\ A\ge 0,B\ge 0,$ such that $\rank A\ge 2,\rank
B\ge 2,$ and let 
$$D:=\{x\in \RR^n:Q_A(x)< 1\},\  E:=\{x\in \RR^n:Q_B(x)< 1\}.
$$
If the boundaries of $D$ and $E$ don't intersect
transversally anywhere, and if $E\cap D\not=\emptyset$, then either 
$E\subset D$ or $D\subset E$.
\end{proposition}

\noi {\bf Proof.} We can argue similarly as in the proof of the case
$n\ge 3$ of Proposition \ref{8M}. 

Let $\R_A=\ker A$ and $\R_B=\ker B$ denote the radical of $Q_A$ and
$Q_B.$ Since $\rank A\ge 2,$ we can find vectors $\xi_1,\xi_2\in\RR^n$
such that $\xi_1\wedge \xi_2\ne 0$ and $\R_A\subset
\xi_1^\perp\cap\xi_2^\perp,$ and similarly $\tilde
\xi_1,\tilde\xi_2\in\RR^n$ such that $\tilde\xi_1\wedge \tilde\xi_2\ne
0$ and $\R_B\subset \tilde\xi_1^\perp\cap\tilde\xi_2^\perp.$ 

If neither $E\subset D$ nor $D\subset E$, we may choose points
$x_0=0\in E\cap D$, $x_1\in E\setminus D$ and $x_2\in D\setminus 
E$ spanning a two-dimensional  plane $V.$ Choose linearly independent
vectors $\eta_1,\dots,\eta_{n-2}\in\RR^n$ such that
$V=\eta_1^\perp\cap\cdots\cap\eta_{n-2}^\perp,$ and consider
$\om:=\xi_1\wedge\xi_2\wedge\eta_1\wedge\cdots\wedge\eta_{n-2}.$ If
$\om\ne 0,$ then $V\cap\R_A=\{0\},$ so that $Q_A|_V>0$ is a positive
definite, i.e., elliptic, quadratic form on $V.$ Similarly, if also 
$\tilde\om:=\tilde\xi_1\wedge\tilde\xi_2\wedge\eta_1\wedge\cdots
\wedge\eta_{n-2}\ne 0,$ then $Q_B|_V>0$ is elliptic, and we can
conclude the proof as in Proposition\ref{8M}.

We shall show that we can slightly vary the points $x_1$ and $x_2$  in
order to achieve that  $\om\ne 0$ and $\tilde\om\ne 0.$
Indeed, since $E\setminus D\ne \emptyset,$ and since the boundary of
$D$ is a smooth submanifold of codimension 1, $E\setminus D$ has  
non-empty interior, and we may assume that $x_1$ lies in this
interior. Similarly, we may assume that $x_2$ lies in the interior of
$D\setminus E.$ Therefore, there is some $\ve>0,$ such that
$x_1+y_1\in E\setminus D$ and $x_2+y_2\in D\setminus E $ for all
$y_1,y_2\in B_\ve(0).$ Let us now replace $\eta_j$ by
$\tilde\eta_j=\eta_j+\delta_j,$  with $\delta_j\in\RR^n$ sufficiently
small, in order to achieve that $\om\ne 0,\tilde\om\ne 0,$ for the
corresponding $n$-forms $\om$ and $\tilde\om.$ Then we can find
$y_1,y_2\in B_\ve(0)$ such that
$\tilde\eta_j(x_1+y_1)=0=\tilde\eta_j(x_2+y_2)$ for $j=1.\dots,n-2$ 
(just solve the linear equations $\tilde\eta_j\cdot y_1=-\delta_j\cdot
x_1,$ and  $\tilde\eta_j\cdot y_2=-\delta_j\cdot x_2,$ for
$j=1,\dots,n-2.$) Replacing $x_1$ by $x_1+y_1$ and $x_2$ by $x_2+y_2,$
we can thus assume that $\om\ne 0$ and $\tilde\om\ne 0,$ which
completes the proof.

\qed
\medskip 

Let us next extend the previous results to non-semidefinite forms. 

 If $W$ is a linear subspace of $\Sym(n,\RR),$ we say that ${A}\in W$
has {\it maximal rank in W}, if  $\rank A=\maxrank M.$  

For  $A,B\in \Sym(n,\RR)$  and $r\ge 0,$   let 
\[
\Gamma^r_{{A}}:=\{Q_{{A}}\le
r\},\quad\Gamma^r_{{B}}:=\{Q_{{B}}\le r\}\, .
\]
Then in particular $\Gamma_{A}=\Gamma_{A}^0,
\Gamma_{B}=\Gamma_{B}^0.$ Notice that if  $z\in\de\Gamma_{{A}}^r, $ and
if ${A} z\ne0,$ then $\de\Gamma_{{A}}^r$ is an analytic manifold of
codimension $1$ near $z,$ and ${A} z$ is normal to it at $z.$  If such a
vector $z$ exists, we say that $\de\Gamma_{{A}}^r$ is a variety of
dimension $n-1.$ If $z\in\de\Gamma^r_{{A}}\cap\de\Gamma^r_{{B}},$  and
if ${A} z$ and ${B} z$ are linearly independent, then we say that the
boundaries $\de\Gamma^r_{{A}}$ and $\de\Gamma^r_{{B}}$ {\it intersect
transversally} at $z.$

\begin{theorem}\label{8N}
Let $W\subset \Sym(n,\RR)$ be a non-trivial subspace of
$\Sym(n,\RR)$ and $r\ge 0,$  and fix 
$A\in W $ of maximal rank in $W,$ with signature
$(k,\ell_1).$ Let also   $B\in W,$  and assume  that
$\de\Gamma_{{A}}^r$ and
$\de\Gamma_{{B}}^r$ are varieties of dimension $n-1.$

If there is no point in
$\de\Gamma^r_{{A}}\cap\de\Gamma^r_{{B}}$ at which the
boundaries $\de\Gamma^r_{{A}}$ and $\de\Gamma^r_{{B}}$ intersect
transversally, and if $\parallel {B}-{A}\parallel$ is sufficiently
small, then either $\Gamma^r_{{A}}\subset\Gamma^r_{{B}}$ or
$\Gamma^r_{{B}}\subset\Gamma^r_{{A}}$, provided $k\ge 2.$
\end{theorem}

\noi {\bf Proof.} Let $\ell:=n-k.$ After rotating the coordinates and
scaling in every coordinate, if necessary, we may assume without loss
of generality that
\[
{A}=\left(\begin{array}{cc}I_k & 0\\
                          0   & -A_2\end{array}\right)\,
                          ,\quad
{B}=\left(\begin{array}{cc}B_1 & B_3\\
                           \trans B_3 & -B_2 \end{array}\right)
\]
with respect to the decomposition $\RR^n=\RR^k\times\RR^{\ell}$,
where, by our assumption, $k\ge 2,$ and where $A_2\ge 0$ has rank
$\ell_1.$\par%%

If $\ell=0$, i.e., if ${A}=I_n$, and if ${B}$ is so close to ${A}$
that ${B}=B_1>0$, then the statement is clear, by Proposition \ref{8M},
if $r>0$, and the case $r=0$ is trivial.\par%%

\medskip

So, assume that $\ell\ge 1$.
We also assume that ${B}$ is so close to ${A}$ that $B_1>0.$ 

As in the proof of Theorem \ref{8C}, we decompose $z=(x,y)$,
$x\in\RR^k$, $y\in\RR^{\ell}$, and put, for $y\in\RR^{\ell}$,
\begin{eqnarray*}
Q^y_{{A}}(x) & := & Q_{{A}}(x,y) = |x|^2-\al(y)\, ,\\
Q^y_{{B}}(x) & := & Q_{{B}}(x,y) =  \trans xB_1x+2(B_3y)\cdot
x-\beta(y),\end{eqnarray*}
where
$\al(y)  :=  \trans y A_2 y\ge 0$ and $\beta(y)  :=   \trans y B_2 y.$
Notice that  we can write 
\begin{equation}\label{qe}
Q_B^y(x)=\trans(x+B_4y)B_1(x+B_4 y)-\trans y\tilde B_2 y,
\end{equation}
with $B_4:=B_1^{-1}B_3$ and $\tilde B_2:=B_2+\trans B_3 B_1^{-1} B_3.$
Observe that $\rank\tilde B_2\le \ell_1,$  since ${A}$ has
maximal rank in $\span\{{A},{B}\}.$ 
\smallskip

 We claim that $\tilde B_2\ge 0$ and $\rank \tilde B_2=\ell_1,$ if
$||{B}-{A}||$ is assumed sufficiently small. 

\noi Indeed, if $||{B}-{A}||$ is sufficiently small, then $Q_{B_2}|_W>0$
on a space $W$  of dimension $\ell_1,$ and thus
also
$Q_{\tilde B_2}|_W>0.$  If we decompose $\RR^\ell=W\oplus H,$ then in
suitable blocks of coordinates $(w,h)$ subordinate to this
decomposition, we may write
 $\tilde B_2=\left(\begin{array}{cc}D & E\\
                           \trans E & F \end{array}\right), $
where $D>0$ has rank $\ell_1.$ Then $Q_{\tilde B_2}(w,h)=\trans
(w+Sh)D(w+Sh) -\trans h\tilde F h, $ for suitable matrices $S,\tilde
F,$ and since
$\rank
\tilde B_2\le \ell_1,$ necessarily $\tilde F=0.$ This implies the
claim.

\smallskip

%Observe also the fact (which will be used at the end of the proof)
%that, e.g., $\Gamma^r_{{A}}\subset\Gamma^r_{{B}}$ holds for one $r>0$
%if  and only if  ${\Gamma^{r'}_{{A}}}\subset{\Gamma^{r'}_{{B}}}$ for
%every
%$r'>0$; this follows just by scaling.\par%%

Let $\Om:=\RR^\ell,$ if $r>0,$ and $\Om:=\RR^\ell\setminus(\ker
A_2\cup\ker\tilde B_2),$ if $r=0.$ Notice that if $y\in\Om,$ then 
$r+\al(y)>0$ and $r+\tilde \beta(y)  :=  r+ \trans y \tilde B_2 y>0,$ 
so that  the ball
\[
D^y_r := \{Q^y_{{A}}\le r\} = \{x\in\RR^k: |x|^2\le r+\al(y)\}
\]
and the ellipsoid
\[
E^y_r := \{Q^y_{{B}}\le r\}
\]
have non-empty interiors. 

Notice that we can exclude the case where $\ell_1=0$ and $r=0,$ for
then $A_2=0=\tilde B_2,$ so that $\Gamma_{A}^0=\{0\}\times \RR^\ell$
and  $\Gamma_{B}^0=\{(x,y)\in\RR^k\times\RR^\ell:x=-B_4 y\}$ have both
codimension greater or equal to two, so that there is nothing to
prove. 

We therefore assume henceforth that $\ell_1\ge 1,$ or $r>0.$ Then
$\Om$ is dense in $\RR^\ell.$ 

Fix $y_0\in\RR^\ell$ such that $r+\al(y_0)>0.$ If we choose
$||{B}-{A}||$ sufficiently small, then we may assume that
$r+\trans y_0(A_2-\trans B_3 B_1^{-2}B_3)y_0>0.$ Modifying $y_0$
slightly, if necessary, we may in addition assume that $r+\trans
y_0\tilde B_2 y_0>0.$ But then the point $(-B_4 y_0,y_0)$ lies in the
open interior of $D_r^{y_0}$ and of $E_r^{y_0}.$ This shows that the
set 
$$\Om_1:=\{y\in\Om: (D_r^y)^0\cap(E_r^y)^0\ne\emptyset\}$$
is a non-empty, open subset of $\Om.$
\medskip

\noi{\bf 1. Case: $\ell_1\ge 2,$ or $r>0.$}
\medskip

\noi Then $\Om$ is connected. We  prove that then $\Om_1=\Om.$ 

For, otherwise, we find $y\in\Om_1$ and $y'\in\Om\setminus\Om_1.$ 
Connect $y$ and $y'$ in $\Om$ by a continuous path
$\gamma:[0,1]\to\Om,$ and put
\[
D_t:=D^{\gamma (t)}_r\, ,\quad E_t:=E^{\gamma (t)}_r\, ,\quad
t\in[0,1]\, ,
\]
and 
\[
\tau := \inf \{t\in [0,1] : \ga(t)\in
\Om\setminus\Om_1\}\, .
\]
Then $\ga(\tau)\in\Om\setminus\Om_1,$ since $\Om\setminus\Om_1$ is
closed in $\Om.$ Now, observe that for $y\in\Om_1,$ the ball $D_y^r$
and the ellipsoid $E_y^r$  don't intersect transversally at common
points of their boundaries. By Proposition \ref{8M}, since $k\ge 2,$
we then have  
\begin{equation}\label{alt}
D^y_r\subset E^y_r,\, \quad\mathrm{or} \quad E^y_r\subset D^y_r,\quad 
\text{for every} \ y\in\Om_1.
\end{equation}
Thus, for some sequence $\{s_j\}_j $ in $[0,\tau[$ tending to $\tau,$
we have
\begin{equation}\label{alt1}
D_{s_j}\subset E_{s_j}, \quad \text{or}\ E_{s_j}\subset D_{s_j}.
\end{equation} 
By continuity, this implies $(D_\tau)^0\cap(E_\tau)^0\ne\emptyset,$
a contradiction. 

We have thus shown that the alternative \eqref{alt} holds for every
$y\in\Om.$ 

However, if $D^y_r\subset E^y_r$ for every
$y\in \Om$, then
\[
Q_{{A}}(x,y)\le r\quad\Rightarrow\quad Q_{{B}}(x,y)\le r \quad \text{for
all}\  y\in\Om\, .
\]
By continuity, this remains true also for all $y\in\RR^\ell$, so that
in this case $\Gamma^r_{{A}}\subset\Gamma^r_{{B}}$.
Similarly, if $E^y_r\subset D^y_r$ for every $y\in\Om$, then we
find that $\Gamma^r_{{B}}\subset\Gamma^r_{{A}}$.\par

Assume therefore that there exist
$y,y'\in\Om$ such that
\[
D^y_r\subset E^y_r\quad\mathrm{and}\quad E^{y'}_r\subset
D^{y'}_r\, .
\]

We claim that we can then find an
$\eta\in\Om$ such that
\begin{eqnarray}\label{8o}
E^{\eta}_r = D^{\eta}_r.
\end{eqnarray}
Since $\Om$ is connected,  we can connect
$y$ and $y'$ in $\Om$ by a continuous path
$\gamma:[0,1]\to\Om$. Put
\[
D_t:=D^{\gamma (t)}_r\, ,\quad E_t:=E^{\gamma (t)}_r\, ,\quad
t\in[0,1]\, .
\]
Then $D_0\subset E_0$, $E_1\subset D_1$. Let
\[
\tau := \mathrm{inf} \{t\in [0,1] : E_t\subset D_t\}\, .
\]
We claim that $E_{\tau}=D_{\tau}$, which proves \eqref{8o}. 

Since
$E_t$ arizes from a fixed centrally symmetric ellipsoid $E$ by scaling
with a positive factor and translation by some vector, both depending
continuously on $t$ (and similarly for $D_t$), we clearly have
$E_{\tau}\subset D_{\tau}$. The case $\tau=0$ is then obvious, so
assume $\tau>0.$ Then $D_s\subset E_s$ for $s<\tau$, so that
$D_\tau\subset E_\tau, $ by continuity.\par

%But, if $E_{\tau}\subsetneq D_{\tau}$, then $|E_{\tau}|<|D_{\tau}|,$
%where $|M|$ denotes the Lebesgue volume   of a Lebesgue measurable
%subset $M$ of $\RR^k.$ On the other hand,
%$D_s\subset E_s$ for $s<\tau$, so that $|D_s|\le |E_s|$, hence, by
%continuity,
%$|D_{\tau}|\le |E_{\tau}|$, a contradiction.\par

With $\eta$ as in \eqref{8o}, we have in particular
\begin{equation}\label{8p}
\{Q^{\eta}_{{B}}=r\} = \{Q^{\eta}_{{A}}=r\} = \{x:|x|^2 =
r+\al(\eta)\}\, .
\end{equation}
This implies
\[
 \trans xB_1x + 2(B_3\eta)\cdot x-(\beta (\eta)+r)=0\quad \text{for
 every }  x\ \text {satisfying} \  |x| = (r+\al(\eta))^{1/2}\, .
\]
Exploiting this for $x$ and $-x$, we see that
\begin{equation}\label{8q}
B_3\eta = 0\, ,
\end{equation}
and then, by scaling,
\[
 \trans xB_1 x=\frac{\beta(\eta)+r}{\al(\eta)+r}|x|^2\quad \text{for all}\ 
x\in\RR^k\, ,
\]
hence
\begin{equation}\label{8r}
B_1 = \beta_1 I_k,
\end{equation}
with $\beta_1 :=(\beta(\eta)+r)/(\al(\eta)+r)$.\par

Let $z :=(x,\eta),$ where $x= (r+\al(\eta))^{1/2}x',$
with $x'\in S^{k-1}.$ Then $z\in\de\Gamma^r_{{A}}\cap\de\Gamma^r_{{B}}$,
and since ${A} z$ is normal to $\de\Gamma^r_{{A}}$ and ${B} z$ normal to
$\de\Gamma^r_{{B}}$ at $z$, the assumptions in the theorem imply
that ${A} z\wedge {B} z=0,$ where, because of \eqref{8q}, 
$${A} z=\left(\begin{array}{cc}x\\ -A_2\eta
\end{array}\right), \quad 
{B} z=\left(\begin{array}{cc}\beta_1 x\\ 
\trans B_3 x-B_2\eta\end{array}\right).
$$
Thus ${B} z=\beta_1{A} z,$ hence
\begin{equation}\label{16}
 \trans B_3x-B_2\eta+\beta_1 A_2\eta=0\quad \text{whenever}\
|x|=(r+\al(\eta))^{1/2}.
\end{equation}
Taking the scalar product with $\eta,$ in view of \eqref{8q} we get 
$-\beta(\eta)+\beta_1\al(\eta)=0,$ hence $r(\al(\eta)-\beta(\eta))=0.$ 
\medskip

If $r>0,$ this implies $\al(\eta)=\beta(\eta),$ hence 
$$\beta_1=1.$$\par

   If $r=0$, notice that $\beta_1$ is close to 1, in view of
\eqref{8r} and  since $||{A} -{B} ||$ is assumed small. Therefore,  we
may replace ${B}$ by $\beta^{-1}_1{B}$, without loss of generality, so
that again we may assume that  $\beta_1=1.$\par

Moreover, by \eqref{16},  then
$$\trans B_3x=(B_2- A_2)\eta\quad \text{whenever}\
|x|=(r+\al(\eta))^{1/2}.
$$
Applying this to $x$ and $-x,$ we find that 
%If $x,x'\in\RR^k$ have length $R:=(r+\al(\eta))^{1/2},$ then this
%implies $\trans B_3(x-x')=0.$ But, since $k\ge 2,$ such vectors
%$x-x'$ span $\RR^k,$ and thus
 $B_3=0 $ and $A_2\eta=B_2\eta.$  We therefore obtain
\begin{eqnarray}\label{8s}
{A} =\left(\begin{array}{ccc}I_k     & 0  \\
                             0 & -A_2\\
    \end{array}\right), \quad
 {B}= \left(\begin{array}{ccc}I_k     & 0  \\
                             0 & -B_2\\
    \end{array}\right), 
\end{eqnarray}
where then also $B_2\ge 0.$ 

Thus 
$$Q_{A}(x,y)=|x|^2-Q_{A_2}(y),\quad Q_{B}(x,y)=|x|^2-Q_{B_2}(y).$$
Fix $x$ such that $|x|^2-r=1.$ Then 
\begin{eqnarray*}
\{y\in\RR^\ell: Q_{A}(x,y)\le r\}&=&\{y\in\RR^\ell: Q_{A_2}(y)\ge 1\},\\
\{y\in\RR^\ell: Q_{B}(x,y)\le r\}&=&\{y\in\RR^\ell: Q_{B_2}(y)\ge 1\},
\end{eqnarray*}
and the complements of these sets are given by 
$$D:=\{y\in\RR^\ell: Q_{A_2}(y)< 1\},\quad E:=\{y\in\RR^\ell:
Q_{B_2}(y)<1\}.
$$
Since, by our assumptions, the boundaries of $D$ and $E$ don't
intersect transversally, we can apply Proposition \ref{8M0} and find 
that $D\subset E,$ or $E\subset D.$ Assume, e.g., that the first
inclusion holds. Scaling in $y$ and taking complements, we then see
that 
$$\{y\in\RR^\ell: Q_{B_2}(y)\ge s\}\subset \{y\in\RR^\ell:
Q_{A_2}(y)\ge s\}\quad\text{for every}\ s\ge 0,$$
hence
$$Q_{B}(x,y)\le r \implies Q_{A}(x,y)\le r,$$
whenever $|x|^2\ge r.$ This implies $\Gamma_{B}^r\subset\Gamma_{A}^r.$

\medskip

\noi{\bf 2. Case: $\ell_1=1$ and $r=0.$}
\medskip

\noi Then there are vectors $\xi_1,\xi_2\in\RR^\ell\setminus\{0\}$ such
that $A_2=\xi_1\otimes\xi_1$ and $\tilde B_2=\xi_2\otimes\xi_2,$ i.e., 
\begin{eqnarray*}
Q^y_{{A}}(x) & := &|x|^2-(\xi_1\cdot y)^2\, ,\\
Q^y_{{B}}(x) & := & \trans(x+B_4y)B_1(x+B_4
y)-(\xi_2\cdot y)^2,\end{eqnarray*} and $\Om=\RR^\ell\setminus
(\xi_1^\perp\cup \xi_2^\perp)$ is non-connected. 
\medskip

\noi{(a) Assume first that $\xi_1\wedge \xi_2\ne 0.$} 
 Then $\Om$ consists of four connected components, on each of  which
the sign of $\xi_1\cdot y$ and  $\xi_2\cdot y$ is constant. Let $\P$
be one if these components such that $\P$ contains a point in $\Om_1.$ 
Arguing as in Case 1, we can then conclude that $\P\subset \Om_1,$ so
that the alternative \eqref{alt} holds for every $y\in\P.$ In
particular, if $y\in\P$ is sufficiently close to the
hyperplane $\xi_1^\perp,$ then $|D_0^y|<|E_0^y|,$ hence $D_0^y\subset
E_0^y.$ Here,  $|M|$ denotes the Lebesgue volume   of a Lebesgue
measurable subset $M$ of $\RR^k.$ Similarly, if $y'$ is sufficiently
close to the hyperplane
$\xi_2^\perp,$ then  $E_0^{y'}\subset D_0^{y'}.$ Connecting $y$ and
$y'$ by some  continuous path
$\gamma$ in $\P,$ we can thus conclude as in Case 1 that
there is some $\eta\in\P$ such that $E_0^\eta=D_0^\eta.$ As in Case
1, this implies that $B_3=0$ and $B_1=A_1,$
without loss of generality. 
Moreover, by \eqref{16}, since $\beta_1=1,$ then $B_2\eta=A_2\eta, $
hence $(\xi_1\cdot\eta)\xi_1=(\xi_2\cdot\eta)\xi_2.$ This contradicts
our assumption that $\xi_1\wedge \xi_2\ne 0,$ and thus this case
cannot arize. 

\medskip

\noi{(b) Assume finally that $\xi_1\wedge \xi_2= 0.$} By a linear
change of coordinates in $\RR^\ell,$ we may then assume that
$\xi_1\cdot y=y_1$ and $\xi_2\cdot y=ay_1,$ for some $a>0,$ so that
$\Om$ consists of the two connected components
$\Om_{\pm}:=\{y\in\RR^\ell:\pm y_1>0\}.$ As before, at least one of
these components must belong to $\Om_1.$ But, since $D_0^{-y}=-D_0^y$ 
and $E_0^{-y}=-E_0^y,$ we see that $y\in\Om_1$ if and only if
$-y\in\Om_1,$ and thus $\Om=\Om_1.$ In particular, the alternative
\eqref{alt} holds for every $y\in\Om.$ Following the arguments applied
in Case 1, assume again that there are  $y,y'\in\Om$ such that
\[
D^y_r\subset E^y_r\quad\mathrm{and}\quad E^{y'}_r\subset
D^{y'}_r\, .
\]
We can then assume that $y$ and $y'$ ly in the same component of
$\Om,$ say, e.g., $\Om_+.$ 

To see this, notice  that a change of sign of $y$ does not
change the volumes of $D^y_r$ and $E^y_r$. Thus, if $D^y_r\subset
E^y_r$, then $E^{-y}_r\subset D^{-y}_r = D^y_r$ would imply
$E^{-y}_r\subset D^y_r\subset E^y_r$, hence $E^{-y}_r = D^y_r =
D^{-y}_r$. This shows that $D^y_r\subset E^y_r$ implies
$D^{-y}_r\subset E^{-y}_r$.

Then we can connect $y$ and $y'$ within $\Om_+$ by some continuous path
$\gamma, $ and thus  find  some $\eta\in\Om_+$ such that
$E_0^\eta=D_0^\eta.$ As before, this implies that $B_3=0$ and
$B_1=A_1,$without loss of generality, and moreover $B_2\eta=A_2\eta,
$ hence  $a=1.$ But then $Q_{A}=Q_{B},$ hence $\Gamma_{A}^0=\Gamma_{B}^0.$

\qed

\begin{remark}\label{8counter}
{\rm If $k=1,$ the statement in Theorem \ref{8N} may fail to be true.
Take, for instance, $Q_{A}(x,y):= x^2-y^2,\ Q_{B}(x,y):= x^2-2\ve xy -y^2,
\ \ve>0,$ for $(x,y)\in\RR^2.$}
\end{remark}

\begin{cor}\label{8t}
Let $W\subset \Sym(n,\RR)$ be a non-trivial subspace of
$\Sym(n,\RR),$   and fix 
$A\in W $ of maximal rank $r$ in $W.$ Let also   $B\in W,$  and
assume  that
$\de\Gamma_{{A}}$ and
$\de\Gamma_{{B}}$ are varieties of dimension $n-1.$

%Let ${A},{B}\in\Sym(n,\RR),$ and assume that ${A}$ has
%maximal rank $r$ in $\span\{{A},{B}\}.$ 
%Assume also that $\de\Gamma_{{A}}$ and $\de\Gamma_{{B}}$ are
%varieties of dimension $n-1.$
 
If there is no point in
$\de\Gamma_{{A}}\cap\de\Gamma_{{B}}$ at which the
boundaries $\de\Gamma_{{A}}$ and $\de\Gamma_{{B}}$ intersect
transversally, and if $||{B} -{A} ||$ is sufficiently small, then
either $\Gamma_{{A}}\subset\Gamma_{{B}}$ or
$\Gamma_{{B}}\subset\Gamma_{A}$, provided $r\ge 3$.
\end{cor}

\noi {\bf Proof.} Assume that ${A}$ has signature $(k,\ell_1)$. The case
$k\ge 2$ is then covered by Theorem \ref{8N}. 

If $k\le 1$, then
$\ell_1\ge 2$, since $r=k+\ell_1\ge 3$. The case $k=0$ is trivial,
since then
$\Gamma_{{A}}=\Gamma_{{B}}=\RR^n$. So assume $k=1$. Applying Theorem
\ref{8N} to $-{A}$ and $-{B}$, we find that
$\Gamma_{-{A}}\subset\Gamma_{-{B}}$ or
$\Gamma_{-{B}}\subset\Gamma_{-{A}}$. Say, the first inclusion holds.
Then also $(\Gamma_{-{A}})^0\subset(\Gamma_{-{B}})^0$, hence
\[
\Gamma_{{B}}=\RR^n\setminus (\Gamma_{-{B}})^0\subset\RR^n\setminus
(\Gamma_{{A}})^0=\Gamma_{{A}}\, ,
\]
so that the statement of the corollary is true also if $k=1$.\par
\qed

From Theorem \ref{8C} and Corollary \ref{8t} we immediately obtain

\begin{theorem}\label{8U}
Let ${A},{B}\in\Sym(n,\RR)$ be a non-dissipative pair,  and let us assume
that
$\maxrank\{{A},{B}\} \ge 3.$ 

If ${A}$ and ${B}$ are linearly independent, then there exists a point
$z\in\de\Gamma_{{A}}\cap\de\Gamma_{{B}}$ at which the
boundaries $\de\Gamma_{{A}}$ and $\de\Gamma_{{B}}$ intersect
transversally.
\end{theorem}

\noi {\bf Proof.} Let $W:=\span_\RR\{{A},{B}\},$ and choose a basis 
$\tilde{A},\tilde{B}$ of $W$ such that $\rank
\tilde{A}=\maxrank\{{A},{B}\},$  and such that
$||\tilde{B}-\tilde{A}||$ is so small that  Corollary
\ref{8t} applies to the pair $\tilde {A},\tilde {B} $ (notice that the
boundaries $\de\Gamma_{\tilde {A}}$ and
$\de\Gamma_{\tilde{B}}$ have dimension $n-1,$ because $\tilde {A}$ and
$\tilde {B}$ have signature $(k,\ell_1),$ with $k\ge 1$ and $\ell_1\ge
1,$ since $W$ is non-dissipative.)  Observe that 
$\de\Gamma_{{A}}\cap\de\Gamma_{{B}}=\de\Gamma_{\tilde
{A}}\cap\de\Gamma_{\tilde{B}}.$ 

Assume now that  there is no point $z\in
\de\Gamma_{{A}}\cap\de\Gamma_{{B}}$ at which ${A} z\wedge{B} z\ne 0.$ Then
there is also no point
 $z\in\de\Gamma_{\tilde {A}}\cap\de\Gamma_{\tilde {B}}$ at which $\tilde
{A} z\wedge\tilde {B} z\ne 0,$ and consequently  $\Gamma_{\tilde
{A}}\subset
\Gamma_{\tilde {B}}$ or $\Gamma_{\tilde {B}}\subset
\Gamma_{\tilde {A}}.$ However, in view of Lemma \ref{8A}, after a
suitable linear change of coordinates we may assume that $\tilde {A}$
and $\tilde {B}$ have vanishing traces. Then, by Theorem \ref{8C},  
$\tilde {A}$ and $\tilde {B}$ are linearly dependent, hence so are 
${A}$ and ${B}.$ This proves the theorem.

\qed

\begin{remark}\label{2.10}
{\rm The analogous statement is false for $r=2$. Take,
for example, $Q_{{A}}(x,y)=x^2-y^2, \ Q_{{B}}(x,y)=xy,\ 
(x,y)\in\RR^2.$}
\end{remark}

\setcounter{equation}{0}
\section{The form problem}

We begin with some auxiliary results and background information on
semi-algebraic sets, which will be useful later. 
 
\subsection{Auxiliary results}

\begin{lemma}\label{8X} Let $A,B\in\Sym (n,\RR)$ such that $A$ is
not semi-definite, and assume that
$\de\Gamma_A\subset\de\Gamma_B$. \be \item[(a)] Then there is a
constant $c\in\RR$ such that $B=cA$. \item[(b)] If
$\mathrm{rank}\, A=2$, so that $\de\Gamma_A\setminus\{0\}$ is not
connected, then let us choose linear coordinates $x=(x_1,\ldots,
x_n)$ so that $Q_A=ax_1x_2$, with $a\in\RR\setminus\{0\}$. If
$Q_B$ vanishes on one of the components of
$\de\Gamma_A\setminus\{0\}$, say $Q_B(x)=0$, if $x_1=0$, then
there is some $b\in\RR^n$ such that $Q_B(x)=(b\cdot x)x_1$. \ee
\end{lemma}

\noi {\bf Proof.} (a) After applying a suitable linear change of
coordinates, we may assume that we can split coordinates
$x=(u,v,w)\in\RR^k\times\RR^{\ell}\times\RR^m=\RR^n$ such that
\[
Q_A(x)=|u|^2-|v|^2=:Q_{A_1}(u,v)\, .
\]
Setting $y:=(u,v)\in\RR^{k+\ell}$, we can then write $Q_B$ as
\[
Q_B (y,w)=Q_{B_1}(y)+(B_2y)\cdot w+Q_{B_3}(w)\, ,
\]
with $B_1\in\Sym (k+\ell,\RR)$, $B_3\in\sgn (m,\RR)$ and $B_2$ a
real $m\times (k+\ell)$-matrix.
Let $y\in\de\Gamma_{A_1}$. Then $Q_B (y,w)=0\; \text{for all}\ \;
w\in\RR^m$, hence
\begin{equation}\label{8y}
B_3=0,\  B_2 y=0\quad\mathrm{and}\quad Q_{B_1} (y)=0\, .
\end{equation}
Since, by our assumptions, $k\ge 1$, $\ell\ge 1$,
$\de\Gamma_{A_1}\,\mathrm{spans}\,\RR^{k+\ell}$. To see this,
choose unit vectors $e_1,\ldots ,e_k\in\RR^k$ spanning
$\RR^k$ and unit vectors  $f_1,\ldots ,f_\ell\in\RR^\ell$ spanning
$\RR^\ell. $ Then the vectors $e_i\pm f_j$ ly in $\de\Gamma_{A_1}$ and
$\mathrm{span}\,\RR^{k+\ell}$. Then \eqref{8y} implies $B_2=0$, so
that
\[
Q_B (y,w)=Q_{B_1}(y)\, .
\]
We may thus reduce ourselves to the case $m=0$. Let us then write
(with new matrices $B_j$)
\[
Q_B(u,v)=Q_{B_1}(u)+Q_{B_2}(v) + 2(B_3 v)\cdot u\, .
\]
For $(u',v')\in S^{k-1}\times S^{\ell-1}$ and $s,t\in\RR$, we have
\begin{eqnarray*}
q (s,t) &:=&
Q_B(su',tv')=s^2Q_{B_1}(u')+t^2Q_{B_2}(v')+2st(B_3v')\cdot u'\\
        &=& \al s^2+\beta t^2+2\gamma st.
\end{eqnarray*}
By our assumptions, $q(s,t)=0$, if $|s|=|t|$. In particular,
$q(t,t)=q(t,-t)=0$, so that $\gamma=0$ and $\al +\beta =0$. Thus
\[
(B_3 v')\cdot u'=0 \quad \text{for all}\ \, u'\in S^{k-1},\, v'\in
S^{\ell-1}\, ,
\]
so that $B_3=0$. Moreover,
\[
Q_{B_1}(u')=-Q_{B_2}(v')\quad \text{for all}\ \, u'\in S^{k-1}, v'\in
S^{\ell-1}\, .
\]
Scaling, this implies
\[
Q_{B_1}(u)=-|u|^2 Q_{B_2}(v')\quad \text{for all}\ \, u\in\RR^k\, ,
\]
hence $B_1=cI_k$, for some $c\in\RR$. Then $Q_{B_2}(v')=-c$, hence
$Q_{B_2}(v)=-c|v|^2\; \text{for all}\ \, v\in\RR^{\ell}$, so that
$B_2=-cI_{\ell}$ and $Q_B=cQ_A$.
\medskip

(b) If $Q_B(x)=0$ for $x_1=0$, then put
$s:=x_1$, $v:=(x_2 ,\ldots ,x_n)$, and write
\[
Q_B(x)=Q_{B_1}(v) +s\beta\cdot v+\gamma s^2\, .
\]
Since $Q_B (v,s)=0$, if $s=0$, we have $B_1=0$, so that
$Q_B=s(\beta\cdot v+\gamma s)$.%%

\qed

\bigskip

%The following lemma will be useful in the sequel.

If $V$ and $W$ are finite dimensional $\KK-$vector spaces, we shall
denote by $L(V,W)$ the space of all linear mappings from $V$ to $W.$
If $V=\RR^k$ and $W=\RR^{\ell}$ are Euclidean spaces, we shall
often identify  a linear mapping $T\in L(V,W)$ with the
corresponding $n\times k-$ matrix in $M^{n\times k}(\KK)$ with
respect to the canonical bases of these spaces, without further
mentioning.  

\begin{lemma}\label{9'}
Let $A\in L(\RR^n,\RR^n)$ and $T\in L(\RR^{n-m},\RR^n)$ be linear
mappings, and assume that $T$ is injective. Then
$$ \rank (\trans TAT)\ge \rank A -2m.
$$
\end{lemma}

\noi {\bf Proof.} It suffices to prove that 
\begin{equation}\label{r1}
\dim(\ker \trans T\, \trans A \,T)\le \dim(\ker \trans A)+m,
\end{equation}
for then $\rank (\trans TAT)=n-m-\dim(\ker \trans T\,\trans A \,T)\ge 
n-\dim(\ker \trans A)-2m=\rank A -2m.$ 

Write $B:=\trans A.$ Since $T$ is injective, \eqref{r1} will follow
from 
\begin{equation}\label{r2}
\dim (\ker \trans T B)\le \dim (\ker B) +m.
\end{equation}
Let $K:=\ker \trans T.$ Then $y\in\ker \trans T B$ if and only if
$By\in K,$ i.e., $\ker \trans T B=B^{-1}(K),$ where clearly $\dim
B^{-1}(K) \le \dim (\ker B) +\dim K=\dim (\ker B) +m.$ This gives
\eqref{r2}.

\qed

\begin{lemma}\label{dot}
Let $A, B \in L(\RR^n, \RR^n)$ be linear symmetric mappings,
i.e., $\trans A=A$ and $\trans B=B,$ and let $I$ be a non-empty
open interval in $\RR$ and 
$E:I\to L(\RR^{n-m},\RR^n)$  a differentiable  mapping such
that
$E(t)$ is injective for every $t\in I.$ Assume  that 
\begin{equation}\label{dot1}
\trans E(t) B E(t)=\mu(t) \trans E(t) A E(t) \quad
\text{for every}\  t\inÊI,
\end{equation} 
where $\mu: I\to \RR$ is a differentiable mapping. Then
$\mu $ is constant, provided $\rank A> 4m.$  

\end{lemma}

\noi {\bf Proof.} Let  $r:=\rank A,$ and put 
$$A(t):=\trans E(t) A E(t), \ B(t):=\trans E(t) B E(t),\quad 
t\inÊI.
$$ 
Then
\[
B(t)=\mu(t) A(t)\quad \text{for all}\  t\in I\, .
\]

Assume that $\mu$ is non-constant. Then there is some point $t_0\in
I$ such that $\dot\mu(t_0):=\frac {d\mu}{dt}(t_0)\ne 0.$
Translating coordinates in $\RR,$ if necessary, we may assume that
$t_0=0.$ Moreover,  replacing then 
$B$ by
$B-\mu(0)A$, we may assume that 
$\mu(0)=0$, hence
$B(0)=0.$ Then
\[
\dot B(0)=\dot\mu(0)A(0)\, ,
\]
hence, by Lemma \ref {9'}, 
\begin{equation}\label{8vv}
\rank\, \dot B(0)=\rank A(0)\ge r-2m\, .
\end{equation}
On the other hand, 
\[
\dot B(0)= \trans E(0) B\dot E(0)+ \trans ( \trans E(0)B\dot E(0))\, ,
\]
and
\[
0=B(0)= \trans E(0)BE(0)\, .
\]
Since $E(0)$ is injective, this implies $ \trans
E(0)B|_{V(0)}=0,$ where $V(0)$ denotes the $(n-m)-$ dimensional
range of $E(0),$  so that  $\rank  \trans E(0)B\le m$, hence
\[
\rank  \trans E(0)B\dot E(0)\le m\, .
\]
The same estimate holds for $ \trans ( \trans E(0)B\dot E(0))$, and
thus
$$
\rank \dot B(0)\le 2m\, .
$$
Together with \eqref{8vv}, this yields $r\le 4m$. 

\qed

\bigskip

In the sequel, we shall consider three matrices 
$A,B,C\in\Sym(n,\RR),$ and shall work under the following assumptions,
unless stated explicitly otherwise:

\begin{assumption}\label{ass}
 $A$ and $B$ are linearly independent and form a
non-dissipative pair, and 
$r:=\maxrank \{A,B\}\ge 3$. Moreover,
 %the joint radical $\ker A\cap \ker B$ of $Q_A$ and $Q_B$ is 
%trivial, i.e., $\ker A\cap \ker B=\{0\},$  
$C$ statisfies the following property:  
\begin{equation}\label{H1}
 Q_C\ \text{vanishes on} \  \{Q_A=0\}\cap\{Q_B=0\}\, .\tag{H1}
\end{equation} 
Unless stated otherwise, we also assume that
$ \rank A=r,$ and that 
$||B-A||$ is sufficiently small.

\end{assumption}

We shall be mostly interested in the situation where
$\RR^n=\RR^{2d}$ is the canonical symplectic vector space, and
$Q_C$ is the Poisson bracket
\begin{equation}\label{H2}
 Q_C=\{Q_A,Q_B\} \tag{H2}
\end{equation} 
of $Q_A$ and $Q_B$.
Our aim will be to prove that, under these assumptions, 
\begin{equation}\label{8v}
C=\al A+\beta B
\end{equation}
for some $\al ,\beta\in\RR,$ provided $r$ is sufficiently large.\\%%

Let us first have a closer look at the structure of the sets 
$\{Q_A=0\}=\de\Gamma_A,$ $\{Q_B=0\}=\de\Gamma_B$ and their intersection 
$$\V :=\de\Gamma_A\cap\de\Gamma_B.$$
 As a reference for the following
results, we recommend \cite{benedetti}.\\

The sets $\de\Gamma_A,\de\Gamma_B$ and $\V$ are real-algebraic,
hence semi-algebraic, so that they admit a finite stratification
into connected, locally closed subsets which are  real-analytic 
subma\-nifolds, each of which is a semi-algebraic set. Recall that the
{\sl dimension} of a semi-algebraic set is the maximal dimension of its
strata.\par%%
If we assume that $B$ is so close to $A$ that also $\rank B=r,$ then
clearly $\de\Gamma_{A}$ and $\de\Gamma_{B}$ are of dimension $n-1$, so
that $\dim\V\le n-1$.\\

Theorem \ref{8U} implies that  there exist  points
$z\in\V$, at which the boundaries of $\de\Gamma_A$
and $\de\Gamma_{B}$ intersect transversally.\par%%
Denote by $\N\subset\V$ the set of all those points $z$ in
$\V$. Then $\N$ is a non-empty  analytic submanifold of
codimension $2$ in $\RR^n$, and $Az,  Bz$ span  the
normal space $N_z\N$ to $\N$ at every point $z\in\N$.\par

In particular, $\N$ decomposes into a finite number of strata of
dimension $n-2$ contained in $\V$, so that $\dim\, \V\ge
n-2$.
It is easy to see that $\dim\,\V =n-1$ is not possible.\par

For, if $\dim\,\V =n-1,$ then there exist a non-empty open neighborhood
$U$ of some point $z\in{\V}$ in $\RR^n$ such that
\[
\de\Gamma_A\cap U=\de\Gamma_B\cap U=\de\Gamma_A\cap\de\Gamma_B\cap
U\, .
\]
However, we may apply a linear change of coordinates so that
\[
Q_A(x)=(x_1^2+\cdots +x^2_k)-(x^2_{k+1}+\cdots +x^2_{\ell_1})\, .
\]
If $k\ge 2$ and $\ell_1\ge 2$, then $\de\Gamma_A\setminus\{0\}$ is
connected, and since $Q_B$ vanishes on
$\de\Gamma_A\setminus\{0\}$ within $U$, it vanishes on the whole
of $\de\Gamma_A,$ by analyticity.
And, if, e.g., $k=1$, then the semi-cones 
$\de\Gamma^\pm_A:=\{x:Q_A(x)=0\,\mathrm{and}\,\pm x_1>0\}$ are
connected, and a similar argument as before shows that $Q_B$
vanishes on at least one of the sets $\de\Gamma^+_A$ or
$\de\Gamma^-_A$. Since $Q_B(-x)=Q_B(x)$, then again $Q_B$ vanishes
on $\de\Gamma_A$.\par%%

Reversing the r\^oles of $A$ and $B$, we see that
$\de\Gamma_A=\de\Gamma_B$. But then $\N=\emptyset$, a
contradiction. We thus have%%
\begin{equation}\label{8w}
\dim\, \V=n-2\, .
\end{equation}

\bigskip

\subsection{The case where no stratum of $\N$ spans $\RR^n.$}

We first  rule out the possibility that no component of $\N$ spans
$\RR^n$.

\begin{theorem}\label{8Z}
Let $A,B,C\in \Sym(n,\RR)$ satisfy our Standing Assumptions
\ref{ass}. Assume that $\maxrank\{A,B\}\ge 7,$  that
$\RR^n=\RR^{2d}$ is symplectic, that $Q_C$ satisfies
H\"orman\-der's bracket condition (H2) and that the joint radical
$\R_{A,B}:= \ker A\cap
\ker B$ of $Q_A$ and $Q_B$ is either trivial, i.e., $\ker A\cap
\ker B=\{0\},$  or a symplectic subspace of $\RR^n.$ 

 Suppose also that no connected component of $\N$ spans
$\RR^n.$ Then $C$ is a linear combination of $A$ and $B$.
\end{theorem}

\noi {\bf Proof.} Without loss of generality, we may assume that 
$\R_{A,B}$ is trivial. For, if $\R_{A,B}$ is a symplectic subspace,
then we can choose a complementary symplectic subspace of $\RR^n,$
and reduce everything to this space. 

Let then 
$X_A(x):=2JAx,\, X_B(x):=2JBx$ denote the Hamiltonian vector fields
associated to
$Q_A$ and $Q_B$. Since
\[
X_AQ_B=\{Q_A,Q_B\}=0\quad\mathrm{on}\quad\N\, ,
\]
the field $X_A$ is tangential to $\de\Gamma_B$ at every point of
$\V\subset\N$, and trivially the same holds for $\de\Gamma_A$, so
that by symmetry in $A$ and $B,$
\[
X_A(x),X_B(x)\in T_x\N\quad \text{for all}\ \, x\in\N\, .
\]

Assume now that no connected component of $\N$ spans
$\RR^n$.

\noi Let $\N_0$ be any component of $\N$. Then
\begin{equation}\label{8AA}
\N_0\subset \nu^{\perp}
\end{equation}
for some unit vector $\nu\in\RR^n$. In particular,
$X_A(x),X_B(x)\in\nu^{\perp}$, hence
\[
\nu\in\;\span_\RR\{X_A(x),X_B(x)\}^{\perp}\quad \text{for all}\ \,
x\in\N_0\, .
\]
Here, $\perp$ denotes the orthogonal with respect to the canonical
Euclidean inner product on $\RR^n$.
Then $J\nu\in\span_\RR\{Ax,Bx\}^{\perp}$, hence
\begin{equation}\label{8bb}
J\nu\in T_x\N_0\quad  \text{for every}\ \, x\in\N_0\, .
\end{equation}
Let
\[
V:= (\span_\RR\{\nu ,J\nu\})^{\perp}\, .
\]
Then $V$ is a $J$-invariant, symplectic subspace of $\RR^n$, and
we can choose orthonormal  coordinates $x_1,\ldots
,x_n$, so that
\[
\nu=e_{n-1}, J\nu=e_n\, ,
\]
where $e_1,\ldots ,e_n$ denotes the associated 
basis, and where $e_1,e_3,\dots e_{2d-1},e_2,\dots, e_{2d}$ forms
a symplectic basis of $\RR^n.$ \par%%

Representing $\N_0$ locally as a graph, we see that \eqref{8bb}
implies that $\N_0$ is locally a cylinder with axis $e_n$, over
a basis $\tilde \M\subset\RR^{n-1};$ by analyticity of $Q_A$ and
$Q_B$, we see that this holds globally:
\begin{equation}\label{8cc}
\N_0=\M\times\RR\, ,
\end{equation}
where
\[
\M:=\{y\in\RR^{n-1}:(y,0)\in\N_0\}\subset\RR^{n-1}\, .
\]
Notice also that, since $\N_0\subset\nu^{\perp}=e^{\perp}_{n-1}$,
we have indeed that $\M$ can be considered as an $(n-3)$-
dimensional submanifold
\begin{equation}\label{8dd}
\M\subset\{y\in\RR^{n-1}:y_{n-1}=0\}=:H\, 
\end{equation}
of the $(n-2)$-dimensional hyperplane $H$ in $\RR^{n-1},$ which can
naturally by identified with $V=H\times \{0\}.$ 

Splitting coordinates $x=(y,s),\  y\in\RR^{n-1}, s\in\RR$, we can
write $A$ in the form
\[
A=\left( \begin{array}{cc}
       A_1 & a_2\\ %\hline
    \trans a_2 & a_3 \end{array} \right)\,
,\quad\mathrm{with}\;a_2\in\RR^{n-1}\simeq e_n^{\perp},\,
a_3\in\RR\,.
\]
Then 
$$Q_A (y,s)=Q_{A_1}(y)+2s(a_2\cdot y)+a_3 s^2.$$
Since $Q_A (y,s)=0\ \text{for all}\  s\in\RR$, if $y\in\M$, we find
that
\begin{equation}\label{8ee}
a_3=0,\; a_2\cdot y=0\quad\mathrm{and}\quad
Q_{A_1}(y)=0\quad \text{for all}\  y\in\M\, .
\end{equation}
In particular, $a_2\in\M^{\perp}$.
\medskip

\noi {\bf Case 1}. $\span_\RR \M=H$.\par%%
\medskip

\noi Then $a_2\in H^{\perp}$, hence $a_2=\tau e_{n-1}$ for some
$\tau\in\RR$, so that

\[
A=\left( \begin{array}{c|c}
  & 0\\
     & \vdots\\
     A_1  & 0\\
      & \tau\\ \hline
     0\cdots 0\,\tau    & 0 \end{array} \right)\, .
\]
Interchanging the r\^oles of $A$ and $B$, we see that $B$ is of the
same form
\[
B=\left( \begin{array}{c|c}
  & 0\\
     & \vdots\\
     B_1  & 0\\
      & \sigma\\ \hline
     0\cdots 0\,\sigma    & 0 \end{array} \right)\, .
\] with  $\sigma\in\RR.$ If both $\tau$ and $\sigma $  were zero,
then $e_n$ would be in the joint radical $\R_{A,B}$ of $Q_A$ and
$Q_B,$ which is assumed to be trivial. Thus, at least one of these
numbers is non-zero, and forming suitable new linear combinations
of $A$ and
$B$  (dropping the assumption that $|| B-A ||$ be small),  we may
assume without loss of generality that 

\[
A=\left( \begin{array}{c|c}
                                                 & 0\\
                                                 & \vdots\\
                                            A_1  & 0\\
                                                 & 1\\ \hline
                                0\cdots 0 \,1      & 0 \end{array}
\right)\,
                                , \quad
B=\left( \begin{array}{c|c}
                                             B_1 & 0\\ \hline
                                              0  & 0 \end{array}
                                              \right)\, ,
\]
so that
\begin{eqnarray}\label{8ff}
Q_A(y,s) & = & Q_{A_1}(y)+2sy_{n-1},\\ \nonumber %%
Q_B(y,s) & = & Q_{B_1}(y)\, .
\end{eqnarray}
Since $A$ and $B$ are linearly independent, we have rank $B_1\ge
1$. Moreover, rank $B_1=1$ is not possible, since then $B_1$ would
be semi-definite. Therefore, rank $B_1\ge 2$.
\medskip

\noi {\bf Case 1(a)}. rank $B_1=2$.
\medskip

\noi Since $B_1$ is not semi-definite, we can then find linearly
independent vectors $\eta_1,\eta_2\in\RR^{n-1}$ such that
\[
Q_{B_1}(y)=(\eta_1\cdot y)(\eta_2\cdot y)\, .
\]
Moreover, by \eqref{8ff},
\[  
\V\cap\{(y,s)\in\RR^{n-1}\times\RR:y_{n-1}=0\}=(\de\Gamma_{A_1}\cap\de\Gamma_{B_1}\cap
H)\times\RR\, ,
\]
so that $\M$ must be an $(n-3)$-dimensional stratum of
$\de\Gamma_{A_1}\cap\de\Gamma_{B_1}\cap H$  of codimension 1
in $H$.
In particular, $\M$ is contained in either
$\eta^{\perp}_1$, or in
$\eta^{\perp}_2$. Assume, e.g., that
\[
\M\subset\eta^{\perp}_1\, .
\]

If $\eta_1\wedge e_{n-1}\not= 0$, then $\M$ lies in the subspace
$K:=\{y\in\RR^{n-1}:y_{n-1}=0,\eta_1\cdot y=0\}$ of codimension
$1$ of $H$, hence is an open subset of $K.$ Since $Q_{A_1}$
vanishes on $\M$, we see that $Q_{A_1}$ vanishes on $K$. Write
\[
A_1=\left( \begin{array}{cc}
  A_1' & a_2'\\ %\hline
  \trans a_2  & a_3' \end{array}
  \right)
\]
with respect to the splitting of coordinates $y=(y',y_{n-1})$. Then
$Q_{A_1'}(y')=0$ whenever $\eta_1'\cdot y':=\eta_1\cdot (y',0)=0$.
Lemma \ref{8X} (b) then implies that $Q_{A_1'}(y')=(\eta_1'\cdot
y')(\eta_2'\cdot y)$, for some $\eta_2'\in\RR^{n-2}$. In
particular, rank $A_1'\le 2$. 

On the other hand, if we apply Lemma
\ref{9'}, where $T:\RR^{n-2}\to \RR^n$ is the inclusion mapping $y'\to
(y',0)\in \RR^{n-2}\times \RR^2=\RR^n,$ then we see that $\rank
A_1'\ge \rank A-4=r-4.$ Thus $r\le 6$, in contrast to our
assumptions.\par
\smallskip

Let us therefore assume that $\eta_1\wedge e_{n-1}=0$, say,
without loss of generality, $\eta_1=e_{n-1}$. Then
\[
Q_B(y,s)=y_{n-1}(\eta\cdot y)\, ,
\]
for some $\eta\in\RR^{n-1}$. Write $\eta =\sigma e_{n-1}+v$, with 
$v\in \span_\RR\{e_1,\ldots ,e_{n-2}\}=V$. Then $v\not= 0$, since
$Q_B$ is not semi-definite, so that we can consider $v$ as a
 member of a new symplectic basis of $V$. Replacing $\{e_1,\ldots
,e_{n-2}\}$ by this new basis, we may then assume without loss
of generality that $v=e_1$, i.e.,
\begin{eqnarray*}
Q_B (y,s) & = &
\sigma y^2_{n-1}+y_1y_{n-1}=y_{n-1}(\sigma y_{n-1}+y_1),\\
Q_A(y,s) & = & Q_{A_1}(y)+2sy_{n-1}\, ,
\end{eqnarray*}
and  that the Poisson bracket of $y_1$ and $y_2$ satisfies
$\{y_1,y_2\}=1.$ Since  $\{s,y_{n-1}\}=1$,  we then obtain
$$Q_C =\{Q_A,Q_B\}
    =  2y_{n-1}(2\sigma y_{n-1}+y_1)-y_{n-1}\frac{\de}{\de
    y_2}Q_{A_1}(y).
$$
Thus,
\[
Q_C = y_{n-1}(\gamma\cdot y)\, ,
\]
for some $\gamma\in\RR^{n-1}$.\par 
Assume that  $y_{n-1}\ne 0$, and
$\sigma y_{n-1}+y_1= 0.$ Then $Q_B (y,s)=0\; \text{for every}\ 
s\in\RR.$ Moreover, we can choose $s=s(y)$ such that $Q_A(y,s)=0$.
Then, by \eqref{H1}, $Q_C(y,s)=0$. 
Thus,  $\sigma y_{n-1}+y_1=0$ always
implies $\gamma\cdot y=0$, so that $\gamma=c(e_1+\sigma e_{n-1}),$ for
some $c\in\RR$, hence $Q_C(y)=cy_{n-1}(y_1+\sigma
y_{n-1})=cQ_B(y)$. Thus $C$ lies in the span of $A$ and $B$.
\medskip

\noi {\bf Case 1(b)}. rank $B_1\ge 3$.
\medskip

\noi Since $B_1$ is not
semi-definite, every connected component of
$\de\Gamma_{B_1}\setminus R$ then spans $\RR^{n-1},$ where $R$ denotes
the radical of $Q_{B_1}$. This can be seen by a similar argument as
in the proof of Lemma \ref{8X}.
 Let $\P$ be such a
component of $\de\Gamma_{B_1}\setminus R$, and consider
\begin{eqnarray*}
\tilde {\N}_0:=\{(y,s(y))\in \RR^n:y\in\P, \, y_{n-1}>0\ 
\mathrm{and}\ s(y)=-\tfrac 1 2 Q_{A_1}(y)/y_{n-1}\}.
\end{eqnarray*}
Then $\tilde {\N}_0$ is a connected submanifold of
dimension $n-2$ contained in $\V$, hence contained in a stratum
$\tilde {\N}$ of maximal dimension. By our assumption,
$\tilde {\N}_0$ is also contained in a hyperplane. Since
$\P$ spans, we can then find some vector $\eta\in\RR^{n-1}$ such
that
\[
s(y)+\eta\cdot y=0\quad \text{for all}\  y\in\P^+:=\{y\in
\P:y_{n-1}>0\}, 
\]
hence 
$$Q_{A_1}(y)=2y_{n-1} (\eta\cdot y) \quad \text{for every}\
y\in\RR^{n-1}.
$$
 Put
\[
Q_{D_1}(y):=Q_{A_1}(y)-2y_{n-1}(\eta\cdot y),\  y\in\RR^{n-1}\, .
\]
Then $Q_{D_1}$ vanishes on $\P^+$, hence, by analyticity, also on
the stratum $\Omega$ of $\de\Gamma_{B_1}$ containing $\P^+$, as
well as on $-\Omega$. However, $\Omega\cup  (-\Omega)$ is dense
in $\de\Gamma_{B_1}$, so that $Q_{D_1}$ vanishes on
$\de\Gamma_{B_1}.$ By Lemma \ref{8X} (a), there is thus a constant
$c\in\RR$  such that
\[
Q_{D_1}=cQ_{B_1}\, .
\]
Replacing $A$ by $\frac 1 2(A-cB)$ (possibly dropping the assumption
that $\rank A=r$), we see that
\begin{eqnarray}\label{8gg}
Q_A(y,s) & = & y_{n-1}(\eta\cdot y)+sy_{n-1}=y_{n-1}(\eta\cdot
y+s),\\ \nonumber Q_B(y,s) & = & Q_{B_1}(y)\, .
\end{eqnarray}
Similarly as in the previous case, we can now choose symplectic
coordinates in $V$ such that
\[
\eta\cdot y = \sigma y_{n-1}+\rho y_1\, ,\quad \sigma,\rho\in\RR\, ,
\]
so that
\[
Q_A(y,s)=y_{n-1}(\rho y_1+\sigma y_{n-1}+s).
\]
Then
\begin{eqnarray*}
Q_C & = & \{Q_A,Q_B\}=y_{n-1}\frac{\de}{\de
y_{n-1}}Q_{B_1}(y)-\rho y_{n-1}\frac{\de}{\de y_2}Q_{B_1}(y)\\
& = & y_{n-1}(\gamma\cdot y)\, ,
\end{eqnarray*}
for some $\gamma\in\RR^{n-1}$. We may assume $\gamma\not= 0$.\par

Let then $y\in\de\Gamma_{B_1}$ such that $y_{n-1}\ne 0.$  Then we can
choose
$s=s(y)$, such that $(y,s)\in \V$, hence, by
\eqref{H1}, $y_{n-1}(\gamma\cdot y)=0$.\par

This shows that $\{y\in\RR^{n-1}:Q_{B_1}(y)=0$, $y_{n-1}\not= 0\}$
lies in the hyperplane $\gamma^{\perp}$, so that $\de\Gamma_{B_1}$
is contained in the union of two hyperplanes. This is, however,
not possible, since rank $B_1\ge 3$, so that the strata of
dimension $n-1$ in $\de B_1$ are not flat.
\medskip

\noi {\bf Case 2}. span $\M\subsetneqq H$.

\medskip
\noi  Let $W:=\span_\RR \M$.
Since $\M$ has codimension 1 in $H$, we see that $W$ has
codimension 1 too, so that $\M$ is an open subset of $W$. And,
§$Q_A,Q_B$ vanish on $\N_0=\M\times\RR$, hence also on $W\times
\RR$, so that $\M=W$ and $\N_0=W\times\RR$. Since $\N_0$ is a
linear subspace of codimension 2 in $\RR^n$, we may introduce new
orthonormal coordinates $x_1,\ldots ,x_n$ in $\RR^n$ such that
\[
\N_0=\{x\in\RR^n :x_{n-1}=x_n=0\}\, .
\]

Splitting coordinates $x=(z,y)$, with $z:=(x_1,\ldots ,x_{n-2})\in
\RR^{n-2}$, $y:=(x_{n-1},x_n)\in\RR^2$, we can therefore write
\[
A=\begin{pmatrix}0 & A_2\\
                  \trans A_2 & A_1 \end{pmatrix}\, ,\quad
                 B=\begin{pmatrix}0 & B_2\\
                                  \trans B_2 & B_1 \end{pmatrix}
\]
with respect to these blocks of coordinates. Notice that rank
$A_2\le 2$.\par 

Then, by Lemma \ref{9'}, $r=\rank A\le 4,$ in 
contradiction to  our assumptions. So, under the hypotheses of Theorem
\ref{8Z}, this case cannot arize.

\qed

\begin{remarks}\label{8hh}
{\rm 
\noi (a) The statement of Theorem \ref{8Z} fails to be true, if $r$
to small, e.g., if $r=4$.

 Consider, e.g., the following
counterexample in $\RR^4=\RR^2\times\RR^2$, with coordinates
$(x,y), \  x=(x_1,x_2)$, $y=(y_1,y_2)$, from
\cite{mueller-karadzhov}, Corollary  1.4:
\[
Q_A (x,y):=x_1y_2+x_2 y_1,\quad Q_B (x,y):=x_1 y_1-x_2 y_2\, .
\]
The corresponding matrices $A$ and $B$ are non-degenerate.
Moreover, putting $\xi := \trans (x_1,x_2)$, $\eta := \trans (-y_1,y_2)$ and
$J_2=\begin{pmatrix}0 & 1\\ -1 & 0 \end{pmatrix}$, then
\[
Q_A=\xi\cdot (J\eta),\quad Q_B=-\xi\cdot\eta\, ,
\]
and thus $\{Q_A=0\}\cap \{Q_B=0\} =\{(\xi
,\eta)\in\RR^2\times\RR^2:\xi =0$ or $\eta=0\}$. Therefore, none
of the connected components of $\N$ span $\RR^4.$  However,
\[
Q_C=\{Q_A,Q_B\} = 2(x_1y_2-x_2y_1)\, ,
\]
so that $A,B$ and $C$ do satisfy our standing assumptions and are
linearly independent.
\medskip

\noi  (b) Also the  condition on the joint radical of $Q_A$ and
$Q_B$ is important, as the following example shows:

Let $e_1,\dots,e_d, f_1,\dots,
f_d$ be a canonical symplectic
basis of  $\RR^{2d},\ d\ge 2,$ with
associated coordinates $(x,y)=(x_1,\dots,x_d,y_1,\dots,y_d),$  and
put
\begin{eqnarray}\label{count1}
Q_A(x,y)&:=& x_1^2-(y_1^2+\cdots +y_{d-1}^2+x_2^2+\cdots +
x_{d-1}^2),\\ \notag
Q_B(x,y)&:=& x_1 x_d. 
 \end{eqnarray}
Then $A,B$ form a non-dissipative pair, since, after a suitable
scaling in $x_1,$ we may assume that the matrices corresponding to
$A$ and $B$ have vanishing traces.   Moreover, 
$$Q_C:=\{Q_A,Q_B\}=-2 y_1x_d,$$ 
so that and $A, B $ and $C$ are linearly independent. Then $\rank
A=2(d-1),$ and  $\R_{A,B}=\RR f_d$ is isotropic with respect to the
symplectic form on $\RR^{2d}.$ Moreover, $A,B$ and $C$ are linearly
independent.
Nevertheless, $Q_C$ vanishes on $\{Q_A=0\}\cap \{Q_B=0\}\ .$ 
\medskip

\noi (c) The condition that $Q_C:=\{Q_A,Q_B\}$ cannot be
dropped either, as the next  example demonstrates:
Let 
\begin{eqnarray}\label{count2}
Q_A(x,y)&:=& x_1^2-(y_1^2+\cdots +y_{d}^2+x_2^2+\cdots +
x_{d-1}^2),\\ \notag
Q_B(x,y)&:=& x_1 x_d. 
 \end{eqnarray}
Then again $A,B$ form a non-dissipative pair, as can be seen as
before, and $\R_{A,B}=\{0\}.$  Let 
$$Q_C:= y_jx_d,$$ 
for any $j=1,\ldots,d.$
Then $A, B $ and $C$ are linearly independent, and  $\rank
\{A,B\}=2d.$  Nevertheless,  $Q_C$ vanishes on $\V=\{Q_A=0\}\cap
\{Q_B=0\}=\RR
e_d\cup \{(x,y): x_d=0\  \text{and}\  Q_A(x,y)=0\}.$ Notice that
here the strata of $\V$ of maximal dimension ly in the hyperplane
$\{x_d=0\},$ but there is a lower dimensional stratum not lying in
this hyperplane. 
\medskip

If we slightly modify $Q_A,$ by putting 
$$Q_A(x,y):= x_1^2-(y_1^2+\cdots +y_{d}^2+x_2^2+\cdots +
x_{d}^2),$$
the situation remains the same, only that now $\V$ lies completely
in the hyperplane $\{x_d=0\}.$

}
\end{remarks}

\medskip

In all these examples, $\minrank\{A,B\}=2.$ We shall prove later in
Lemma \ref{rank2} that this is necessarily so, if no stratum of $\N$
spans $\RR^n.$

%A comprehensive discussion of the case $d=2$, including the
%question of local solvability of the associated second order
%operators on $\HH_2,$ would have to extend the work in [KM].\\

\subsection{The case where at least one  stratum of $\N$ spans
$\RR^n.$}

%\end{document}

\bigskip

Let us now go back to our standing assumptions, not requiring that
$\RR^n$ is symplectic and $C$ satisfies (H2). However, let us
assume that
\begin{equation}\label{8ii}
\span_\RR\, \N_0 =\RR^n,
\end{equation}
for some connected component $\N_0$ of $\N$. 

Notice that then, if $U_0$ is a non-empty open subset of $\N_0$,
 also $U_0$ spans $\RR^n$, by
analyticity and connectivity of $\N_0.$ 
 For $x_0\in\N$, we
denote by
\[
V_{x_0}:=T_{x_0}\N=(\span_\RR\, \{Ax_0,Bx_0\})^{\perp}
\]
the tangent space to $\N$ at $x_0$.

Moreover, if $D\in\Sym(n,\RR)$,
then
\[
Q_{D,x_0}:=Q_D|_{V_{x_0}}
\]
denotes the restriction of the quadratic form $Q_{D,x_0}$ to $V_{x_0}$.
Similarly, if $x_0,x_1\in\N$, we put
\begin{eqnarray*}
V_{(x_0,x_1)} & := & V_{x_0}\cap V_{x_1}\, ,\\
Q_{D,x_0,x_1} & := & Q_D|_{V_{(x_0,x_1)}}\, .
\end{eqnarray*}

\begin{lemma}\label{8JJ}
Given $x_0\in\N$, there exist $\alpha (x_0),\beta(x_0)\in\RR$ such
that
\begin{equation}\label{8kk}
Q_{C,x_0} = \alpha (x_0)Q_{A,x_0}+\beta(x_0) Q_{B,x_0}
\end{equation}
and
\begin{equation}\label{8ll}
Cx_0 = \alpha(x_0) Ax_0+\beta(x_0) Bx_0\, .
\end{equation}

The $\alpha(x_0),\beta(x_0)$ may not be unique, but can locally on
$\N$ be chosen as real-analytic functions of $x_0$.
\end{lemma}

\noi {\bf Proof}. Fix $x_0\in\N$. Since $\nabla
Q_A(x_0)=2Ax_0$ and $\nabla Q_B(x_0)=2B x_0$ are linearly
independent, possibly after relabeling the coordinates of $\RR^n$,
we may  assume that
\[
\phi (x) := (x_1,\ldots ,x_{n-2}, Q_A(x),Q_B(x)) =: (x',y)
\]
is a local analytic diffeomorphism near $x_0$. In the new
coordinates $(x',y)$, the point $x_0$ corresponds to $(x'_0,0)$,
and $\N$ to $\{y=0\}$. Let $\tilde {f}:=f\circ \phi^{-1}$, if $f$
is a function defined near $x_0$. Since $Q_C$ vanishes on
$\N$, a Taylor-expansion of $\tilde {Q}_C$ near $(x'_0,0)$ shows
that
\[
\tilde {Q}_C(x',y) =
y_1\tilde {\alpha}(x',y)+y_2\tilde {\beta}(x',y)\,
,
\]
for suitable analytic functions
$\tilde {\alpha},\tilde {\beta}$ defined near
$(x'_0,0).$ Thus, near $x_0$,
\begin{equation}\label{8mm}
Q_C (x) = \alpha(x)Q_A(x)+\beta(x)Q_B(x)\, ,
\end{equation}
for some analytic functions $\alpha ,\beta$ defined near $x_0$.
Taking first derivates at $x_0$, we obtain \eqref{8ll}.\par 

Moreover
applying the second derivative to \eqref{8mm}, we obtain
\begin{equation}\label{8nn}
C=\alpha(x_0)A+\beta(x_0)B+M(x_0)\, ,
\end{equation}
where
\begin{eqnarray}
2M(x_0) & := & \alpha'(x_0)\otimes Ax_0+Ax_0\otimes\alpha' (x_0)\\
\nonumber & & +\beta'(x_0)\otimes Bx_0+Bx_0\otimes\beta'(x_0)\,
.\label{8oo}
\end{eqnarray}
This implies \eqref{8kk}. 

\medskip\qed

\medskip

\noi For $x_0,x_1\in\N_0$, let
\[
N_{(x_0, x_1)}:= \span_\RR\, \{Ax_0,Bx_0,Ax_1,Bx_1\}\, ,
\]
so that $V_{(x_0, x_1)}=N^{\perp}_{(x_0,x_1)}$. Then $\dim\,
N_{(x_0, x_1)}\ge 2$. If $\dim\, N_{(x_0, x_1)}=2\quad \text{for all}\ (x_0,x_1)\in\N^2_0$, then $V_{x_1}=V_{x_0}\quad \text{for all}\ 
x_0,x_1\in\N_0$, so that $\N_0$ would not span $\RR^n$. So, if we put 
$$
m :=  \max_{(x_0,x_1)\in\N^2_0}\,\dim N_{(x_0,x_1)}\, ,
$$
then $m\in\{3,4\}$.

Let us call a subset $U$ of an open domain $\Om\subset\RR^m$ a
{\it generic set} or set of {\it generic points} in $\Om$, if
there exists a non-trivial real analytic function
$f:\Om\rightarrow\RR$ such that $U=\{x\in\Om :f(x)\not=0\}$. A
property will hold for generic points in $\Om$, if it holds for
all points of a generic subset. \par Notice that a generic set in
$\Om$ is open and dense in $\Om$. Clearly, the intersection of a
finite number of generic sets in $\Om$ is again generic.\par
Moreover, if $m=k+\ell$, and if $U$ is a generic subset of
$\Om\subset\RR^k\times\RR^{\ell}$, then let
$\Om^1\subset\RR^k,\, \Om^2\subset\RR^{\ell}$ be any domains
such that $\overline{\Om^1}\times\overline{\Om^2}$ is compact in
$\Om$. Then
$U_1:=\{x\in\Om^1:\, U_x$ is generic in $\Om^2\}$ is generic in
$\Om^1.$  Here $U_x$ denotes the $x$-section 
\[
U_x:=\{y\in\Om^2 :(x,y)\in U\}
\]
of $U$.
\noi  Indeed, if $U=\{(x,y)\in\Om :f(x,y)\not= 0\}$, then
$\Om_1\setminus U_1=\{x\in\Om^1:D^{\alpha}_y
f(x,0)=0\quad \text{for all}\ \alpha\in\NN^{\ell}\}.$ Thus, if we put
\[
g(x):=\sum_{\alpha\in\NN^{\ell}}\ve_{\alpha}D^{\alpha}_y
f(x,0)^2\, ,
\]
for a suitable family at coefficients $\ve_{\alpha}>0$ tending to
$0$ sufficiently fast as $|\alpha |\rightarrow\infty$ ,
then $g$ is a non-trivial real analytic function on $\Om^1$, and
\[
U_1=\{x\in\Om^1 : g(x)\not= 0\}\, .
\]
Analogous definitions and results apply for open domains $\Om$ in
a real analytic manifold.
\medskip

If $m=3$, and if $\dim\,N_{(x_0,x_1)}=3$, then three of the vectors
$Ax_0, B{x_0}, A{x_1}, B{x_1}$ are linearly independent,
say, e.g., $A{x_0}\wedge B{x_0}\wedge A{x_1}\ne 0$. But then
this holds for generic pairs $(x_0,x_1)\in\N^2_0,$ by analyticity
and connectivity of $\N_0$, hence $\dim\, V_{(x_0,x_1)}=n-3$, for
generic $(x_0,x_1)$. A similar argument applies to the case $m=4.$
Therefore, the set
$$\G:=\{(x_0,x_1)\in\N^2_0 :\dim\,
V_{(x_0,x_1)}=n-m\}
$$
is generic in $\N^2_0,$  and
$\dim\,V_{(x_0,x_1)}\ge n- m\ \text{ for all}\ 
(x_0,x_1)\in\N^2_0$.

\begin{proposition}\label{8pp} Assume that $A,B$ and $C$ satisfy
our standing assumptions, and that \eqref{8ii} holds. If
$Q_{A,x_0,x_1}$ and $Q_{B,x_0,x_1}$ are linearly independent for
some $(x_0,x_1)\in\G$, then $C$ lies in the linear span of $A$ and
$B$.
\end{proposition}

\noi {\bf Proof}. Choose connected open neighborhoods $U_j$ of $x_j,
\, j=0,1,$ in $\N_0$ such that $U_0\times U_1\subset\G$. For $(y,z)\in
U_0\times U_1$, we have, by \eqref{8kk},
\[
Q_{C,y,z}=\alpha (y) Q_{A,y,z}+\beta (y) Q_{B,y,z}\, ,
\]
as well as
\[Q_{C,y,z}=\alpha (z) Q_{A,y,z}+\beta (z) Q_{B,y,z}\, .
\]
Moreover, $Q_{A,x_0,x_1}\wedge Q_{B,x_0,x_1}\not= 0$. Shrinking
$U_0$ and $U_1$, if necessary, we may therefore assume that also
$Q_{A,y,z}\wedge Q_{B,y,z}\not= 0\ \text{for every}\  (y,z)\in
U_0\times U_1$, so that
\[
\alpha (y)=\alpha (z),\;\beta (y)=\beta (z)\quad \text{for all}\  (y,z)\in
U_0\times U_1 \, .
\]
Putting $\alpha_0:=\alpha (x_1),\,\beta_0:=\beta (x_1),$ we see that
\begin{equation}\label{8qq}
\alpha (y)=\alpha_0,\; \beta (y)=\beta_0\quad \text{for all}\  y\in
U_0\, .
\end{equation}
But, since $U_0$ is a non-empty open subset of $\N_0$, and since
$\N_0$ spans $\RR^n$, then also $U_0$ spans $\RR^n$, by
analyticity and connectivity of $\N_0.$ \eqref{8qq} and \eqref{8ll}
imply
\[
Cy=\alpha_0 Ay+\beta_0 By\quad \text{for all}\  y\in U_0\, ,
\]
and since $U_0$ spans, we obtain
\[
C=\alpha_0 A+\beta_0 B\, .
\]
\qed

\medskip

There remains the case where
\begin{equation}\label{8rr}
Q_{A,x_0,x_1}\wedge Q_{B,x_0,x_1}=0\quad \text{for all}\ 
 (x_0,x_1)\in\G\, .
\end{equation}

\begin{lemma}\label{8SS}
Let $A$ and $B$ satisfy our standing assumptions, and that there
exists a non-empty domain $U\subset\G$ such that $Q_{A,x,y}$ and
$Q_{B,x,y}$ are linearly dependent for every $(x,y)\in U$. Then
there exists a non-trivial linear combination $D=\alpha A+\beta B$
of $A$ and $B$ such that
\[
Q_{D,x,y}=0\quad \text{for every}\  (x,y)\in U\, ,
\]
provided $r\ge 17$.
\end{lemma}

\noi {\bf Proof}. Since $N_{(x,y)}= \span_\RR\, \{Ax,Bx,Ay,By\}$
varies analytically in $(x,y)\in\G$, shrinking $U$, if necessary,
we may assume that there is an orthonormal basis (constructable,
e.g., by the Gram-Schmidt method)
\[
e_1(x,y),\ldots , e_{n-m}(x,y),\quad (x,y)\in U\, ,
\]
of $V_{(x,y)}$, varying analytically in $(x,y)$. Put
\[
E(x,y):=(e_1(x,y)),\ldots , e_{n-m} (x,y)\in M^{n\times(n-m)}(\RR)\,
,
\]
and
\begin{eqnarray*}
A(x,y) & := &   \trans E(x,y)A E(x,y)\, ,\\
B(x,y) & := &  \trans E(x,y)A E(x,y)\, ,\quad (x,y)\in U.
\end{eqnarray*}
Then $A(x,y)$ and $B(x,y)$ are real $(n-m)\times (n-m)$-matrices,
representing the forms $Q_A|_{V_{(x,y)}}$ and $Q_B|_{V_{(x,y)}}$
with respect to the above  basis of $V_{(x,y)}$. By Lemma
\ref{9'},  we have
\begin{equation}\label{8tt}
\rank\, A(x,y)\ge r-2m
\end{equation}
%Assuming that
%$||B-A||$ is sufficiently small, we have similarly $\rank\,
%B(x,y)\ge r-2m$.\par

If $r\ge 9$, we have in particular that $A(x,y)\not= 0.$
Since $Q_{A,x,y}\wedge Q_{B,x,y}=0$,
there exists thus a unique 
$\alpha(x,y)\in\RR$ such that 
\begin{equation}\label{8uu}
B(x,y)=\alpha(x,y)\; A(x,y),\quad (x,y)\in U\, .
\end{equation}
Since $\alpha (x,y)$ is unique, it depends  analytically on 
$(x,y)\in U.$\par

Assume $\alpha$ is non-constant on $U$. Then we can choose a
differentiable  curve $\gamma :I\rightarrow U$, such that
$\mu :=\alpha\circ\gamma:I\rightarrow\RR$ is non-constant. Put
$E(t):=E(\gamma(t))$, and
\[
A(t):= A(\gamma (t)),\  B(t):=B
(\gamma(t)),\quad  t\in I .
\]
Then 
$$\trans E(t) B E(t)=\mu(t) \trans E(t) A E(t) \quad 
\text{for every}\  t\inÊ]-1,1[,
$$
so that, by Lemma \ref{dot}, $r\le 4m\le 16,$  contradicting our
assumptions. 

Consequently, $\al$ is constant, i.e.,  $\alpha\equiv \alpha_0,$ for
some $\alpha_0\in\RR$. Putting $D:=B-\alpha_0A$, we find that
\[
Q_{D,x,y}=0\, ,
\]
first, for every $(x,y)$ in our shrinked domain $U\subset \G$, but
then also for all $(x,y)$ in our original domain $U$. 

\qed

\begin{lemma}\label{8XX}
Assume that $\N_0\,spans\,\RR^n$, and that $r\ge 12$, and let $U_1$
be a non-empty open subset of $\N_0$. We also assume
\eqref{8rr}.\par

Then, for generic $x_0\in U_1$, there exists an $x_1\in U_1$, such
that for every sufficiently small neighborhood $\tilde U_1$ of $x_0$
and
$U_2$ of $x_1$ in $U_1,$ the following hold:
\medskip

$U:=\tilde U_1\times U_2\subset\G$, and the union of the spaces
$V_{(x,y)}, \ (x,y)\in U,$ spans $\RR^n$.
\end{lemma}

\noi {\bf Proof}. Without loss of generality, we may assume
that $||B-A||$ is so small that $\rank B=\rank A=r.$ 

For generic $x_0\in U_1$, there is some $x_1\in
U_1$ such that $(x_0,x_1)\in\G$. Given such a pair $(x_0,x_1)$,
choose open neighborhoods $\tilde U_1$ of $x_0$
and $U_2$ of $x_1$ in $U_1$ such that $U\subset\G.$ 
\medskip

Assume now that
$\bigcup \{V_{(x,y)}:(x,y)\in U\}$ does not span $\RR^n.$ 
Then there is a unit vector $\nu\in\RR^n$ such that $\nu\perp
V_{(x,y)}$, i.e., $\nu\in N_{(x,y)}, \text{for all}\  (x,y)\in
U$.

We can exclude that $\nu\wedge Ax\wedge Bx=0\ \text{for
every }\  x\in \tilde U_1$. For then 
$\nu\in \span\, \{Ax,Bx\}$, hence $V_x\subset
\nu^{\perp}\quad \text{for every}\  x\in \tilde U_1$. But then
$V_x\subset\nu^{\perp}\quad \text{for all}\  x\in\N_0$, by analyticity and
connectivity of $\N_0$. This implies $z\cdot\nu =
\mathrm{const.}\  \text{for all}\  z\in\N_0$, and since
$0\in\overline{\N_0}$, we have $z\cdot\nu =0$.
Thus $\N_0\subset\nu^{\perp},$ in contradiction to our
assumptions.
\smallskip

 Thus, $\nu\wedge Ax\wedge Bx\not=0$ for generic $x\in
\tilde U_1$. Choose $x_0'\in \tilde U_1$ and an open neighborhood
$U_1'$ of
$x_0'$ in $\tilde U_1$, such that
\begin{equation}\label{8yy}
\nu\wedge Ax\wedge Bx\not=0\quad \text{for every}\  x\in U_1'\, ,
\end{equation}
and put $U':=U_1'\times U_2\subset U$.
\medskip

Assume now first that $m=3$, and that, e.g.,
\[
N_{(x,y)}=\,\span_\RR\,\{Ax,Bx,Ay\}\, ,
\]
first, for $(x,y)=(x_0',x_1)$, and then for every $(x,y)\in U'$
(shrinking $U'_1$ and $U_2,$ if necessary). Since $\nu\in N_{\xy}$,
we have $\nu\wedge Ax\wedge Bx\wedge Ay=0$, hence, by \eqref{8yy},
\[
Ay\in\, \span_\RR\, \{\nu ,Ax,Bx\}\quad \text{for all}\ \xy\in U'\,
.
\]
Fixing $x\in U_1'$, we see that $Ay$ lies in the 3-dimensional
subspace $\span_\RR\,\{\nu ,Ax,Bx\},\ \text{for all}\  y\in U_2$,
so that
$\{Ay:y\in U_2\}$ does not span $\RR^n.$ But then $U_2$ does not
span either, for otherwise, $\span_\RR\,\{Ay:y\in U_2\}=\range A$ 
would have dimension $r,$ hence $r\le 3.$ The case where
$N_{(x,y)}=\,\span_\RR\,\{Ax,Bx,By\}$ can be treated in the same 
way. 
\medskip

There remains the case $m=4$. Then
\[
Ax\wedge Bx\wedge Ay\wedge By\not= 0\, ,
\]
and
\[
\nu\in\, \span_\RR\,\{Ax,Bx,Ay,By\}\quad \text{for all}\  \xy\in
U_1'\times U_2\, .
\]
By \eqref{8yy}, this implies that there exists $(\alpha\xy
,\beta\xy)\in\RR^2\setminus\{0\},$ such that
\begin{equation}\label{8zz}
\alpha\xy Ay+\beta\xy By\in\, \span_\RR\,\{\nu
,Ax,Bx\}\quad \text{for all}\ \xy\in U_1'\times U_2\, .
\end{equation}

Moreover, since $r\ge 9$, the proof of Lemma \ref{8SS}
(see \eqref{8uu}) shows that we may assume that
\[
Q_{B,x_0,x_1}=\alpha_0 \,Q_{A,x_0,x_1},
\]
for some $\alpha_0\in\RR$. Put $D:=B-\alpha_0 A$. Then $D\ne 0$,
and $Q_{D,x_0,x_1}=0$. 
%In suitable orthonormal coordinates, $D$ has then the form
%\[
%D=\left(\begin{array}{l|l} O & D_1\\ \hline
%          \trans D_1 & D_0   \end{array}
%\right)\, ,
%\]
%where $D_0\in\sgn (m,\RR)$. 
By Lemma \ref{9'} then implies that $\rank D\le 2m=8.$ 

Fix $x$ in \eqref{8zz}, and put $W_x:=\span_\RR\,\{\nu ,Ax,Bx\}$.
Then
\[
\alpha\xy Ay+\beta\xy By\in W_x\quad \text{for all}\  y\in U_2\, .
\]
Writing $B=\alpha_0 A+D$, we have, by \eqref{8zz},
\[
(\alpha\xy +\alpha_0\beta\xy) Ay\in W_x+\range D.
\]
If $\alpha\xy +\alpha_0\beta\xy\not=0$ for some $y\in U_2$, then
this implies $r= \rank A\le 3+\rank\, D\le 11$, so that this case
cannot occur.
Assume therefore that $\alpha\xy+\alpha_0\beta\xy
=0\quad \text{for all}\ \xy\in U_1'\times U_2$. Then, by
\eqref{8zz},
\[
\beta\xy (-\alpha_0 A+B)y\in W_x\quad \text{for all}\  y\in U_2\, .
\]
But, since $\rank A=r\ge 12,$ \eqref{8zz} implies that
$\beta\xy$ cannot vanish identically on $U_2$, for generic $x\in
U_1'$, so that
\[
Dy=(-\alpha_0 A+B)y\in W_x,
\]
for generic $\xy\in U_1'\times U_2$. Thus, since $U_2'$ spans
$\RR^n,$
$$
\range D\subset\bigcap_{x \in V_1} W_x=:\W,
$$
for some non-empty open subset $V_1$ in $U_1'$. In particular, rank
$D\le 3$.\par

But, if $W_x=W_y\; \text{for all}\  x,y\in V_1$, then, for fixed $x$,
\[
Ay\in W_x\quad \text{for every}\  y\in V_1\, ,
\]
hence $\range A\subset W_x$, since $V_1$ spans $\RR^n$, a
contradiction. Therefore, $\dim \W\le 2$, hence rank $D\le 2$.
Since $D$ cannot be semi-definite, we thus have $\rank D=2\, .$
Then, there are $\xi,\eta\in\RR^n\setminus \{0\}$ such that
$Q_D(x)=(\xi\cdot x)(\eta\cdot x)$. But this implies
\[
\N_0\subset \xi^{\perp},\;\mathrm{or}\;\N_0\subset \eta^{\perp}\,
,
\]
so that $\N_0$ does not span $\RR^n$, in contrast to our
assumptions. 

\qed

\begin{lemma}\label{8YY}
Let $A,B\in\Sym (n,\RR)$  satisfy our standing
assumptions. Assume further that there is a non-empty open subset
$U\subset\G$ such that the union of the spaces $V_{\xy},\xy\in U$,
spans $\RR^n$, and $r\ge 12$. If $D\in\span_\RR\{A,B\},$ then the
condition
\begin{equation}\label{8ab}
Q_{D,x,y}=0\quad \text{for all}\ \xy\in U
\end{equation}
implies $D=0.$
\end{lemma}

\noi {\bf Proof}. Assume that $D$ satisfies \eqref{8ab}, but that
$D\ne 0.$  Then $\span_\RR\{A,D\} =\span_\RR\{A,B\},$ so that we may
assume without loss of generality that $D=B.$ Notice, however, that
we can then still assume that $\rank A=r,$ but no longer that also
$\rank B=r$ (deviating thus slightly from our standing
assumptions),  since
$B$ then may not be close to
$A.$ We shall show that these assumptions lead to the  
contradiction that $B=0.$

As in the
proof of Lemma
\ref{8SS}, we may assume that there is an orthonormal basis
$e_1\xy,\ldots ,e_{n-m}\xy$ of $V_{\xy}$, varying analytically in
$\xy\in U$.
\medskip

\noi {\bf (a)}  We begin with the case $m=4$. Putting 
$$f_1\xy:=Ax, \ f_2\xy:=Bx,\ f_3\xy:=Ay,\ f_4\xy:=By,$$
 then
\[
e_1\xy,\ldots ,e_{n-4}\xy,f_1\xy,\ldots ,f_4\xy
\]
is a basis of $\RR^n$. Consider the mapping
$F:U\times\RR^{n-4}\rightarrow\RR^n$,
\[
F(x,y,z):=\sum^{n-4}_{j=1}z_je_j\xy,\quad \xy\in U,z\in\RR^{n-4}\,
, \]
where we consider $U$ as an analytic submanifold of
$\RR^{n}\times\RR^{n}$ of dimension $2(n-2).$
\medskip

We shall prove that there is some $(x_0,y_0,z_0)$ such that 
\begin{equation}\label{8ac}
\rank\, DF (x_0,y_0,z_0)=n\, .
\end{equation}
This implies that $F$ is a submersion near $(x_0,y_0,z_0)$, so
that $\bigcup_{\xy\in U}V_{\xy}$ contains a non-empty
open subset $\Om$ of $\RR^n.$ Since, by \eqref{8ab}, $Q_B$ vanishes
on $\Om$, then $Q_B\equiv 0$, hence $B=0.$  \par

To prove \eqref{8ac}, notice that
\[
T_{\xy} U=V_x\times V_y\, ,
\]
and
\[
V_x=\{\xi\in\RR^n:\quad\xi\cdot (Ax)=0,\, \xi\cdot(Bx)=0\}\, .
\]
Consider, for $x,y,z$ fixed, the linear mappings 
\begin{eqnarray*}
\psi_i &:=& V_x\times V_y\to \RR\, ,\\
\psi_i(\xi,\eta)&:=& f_i\xy\cdot (D_{\xy}F(x,y,z)(\xi,\eta)),\quad
i=1,\ldots 4\, .
\end{eqnarray*}
We claim that \eqref{8ac} holds for $(x_0,y_0,z_0):=(x,y,z)$, if and
only if the linear mapping
\[
\Psi:=\left(\begin{array}{c} \psi_1\\ \vdots\\
\psi_4\end{array}\right) :V_x\times V_y\rightarrow\RR^4 
\]
has rank 4.

To this end, observe first that $F(x,y,z)$ is linear in $z,$ so
that 
\begin{equation}\label{help}
DF(x,y,z)(\xi,\eta,\zeta)=D_{\xy}F(x,y,z)(\xi
,\eta)+ F(x,y,\zeta).
\end{equation}
Noreover, if $f_1^*\xy ,\ldots ,f_4^*\xy$ denotes the dual
basis of
$N_{\xy}$ with respect  to the basis $f_1\xy ,\ldots ,f_4\xy$,
i.e., if
\[
f_i^{\star}\xy\cdot f_j\xy=\delta_{ij}\, ,
\]
and if $\psi_i (\xi ,\eta)=w_i$, $i=1,\ldots ,4$, then  let us put
\[
v:=\sum^4_{j=1}w_if_i^{\star}\xy\in N_{\xy}\, .
\]
Then
\[
f_j\xy\cdot [v-D_{\xy}F(x,y,z)(\xi ,\eta)]=0,\quad j=1,\ldots ,4\,
,
\]
so that $v-D_{\xy}F(x,y,z)(\xi ,\eta)\in V_{\xy}$, i.e.,
\[
D_{\xy}F(x,y,z)(\xi ,\eta)-v=-\sum^{n-4}_{j=1}\zeta'_j e_j \xy\, ,
\]
for some (unique) $\zeta' =(\zeta_1',\ldots,
\zeta'_{n-4})\in\RR^{n-4}$. This implies, by \eqref{help}, 
\[
D_{(x,y,z)}F(x,y,z)(\xi ,\eta \,, \zeta
+\zeta')=v+\sum^{n-4}_{j=1}\zeta_j e_j\xy\quad \text{for all} \
\zeta\in
\RR^{n-4}, 
\]
showing that $DF(x,y,z)$ is surjective if and only if  $\Psi$ is
surjective.\par

Observe next that
\[
f_i\xy\cdot F(x,y,z)\equiv 0\,, \quad i=1,\ldots ,4\, ,
\]
so that, by the product rule,
\begin{eqnarray*}
\psi_i(\xi ,\eta) & = & -[(D_{\xy}f_i\xy (\xi ,\eta))]\cdot
F(x,y,z)\\
& = & -f_i (\xi ,\eta)\cdot F(x,y,z)\, .
\end{eqnarray*}
This  easily implies
\begin{eqnarray*}
\Psi (\xi ,\eta)=\left(\begin{array}{c}\psi_1 (\xi ,\eta)\\
\vdots\\
\psi_4 (\xi ,\eta)
\end{array} \right) 
=-\left(\begin{array}{cc} \trans (A\cdot F(x,y,z)) & 0\\
 \trans (B\cdot F(x,y,z)) & 0\\
  0   &  \trans (A\cdot F(x,y,z))\\
  0  & \trans (B\cdot F(x,y,z)) \end{array} \right)
{\xi \choose \eta}\, ,
    \end{eqnarray*}
$ (\xi ,\eta)\in V_x\times V_y\, .$

Notice that, for $\xy\in U$ fixed, there is some $z\in\RR^{n-4}$
such that $\psi_1|_{V_{\xy}}\not= 0$. Indeed, otherwise we would
have 
$(Av)\cdot\xi =0\quad \text{for all}\  v\in V_{\xy}$, $\xi\in V_x$, and in
particular $Q_{A,x,y}=0$. But then rank $A\le 2\cdot 4=8$, hence
$r\le 8$, a contradiction. This shows that $\psi_1\not= 0$ and
$\psi_3\not= 0,$ for generic
$z$, so that
\[
2\le \rank \Psi\le 4,\quad\mbox{for generic}\, z\, .
\]

If rank $\Psi =4$ for some $(x,y,z)$, then \eqref{8ac} holds, and
thus again  $B=0$.\par

So, assume rank $\Psi\le 3\quad \text{for every }\  (x,y,z)\in
U\times\RR^{n-4}$. Then, either the first two rows of $\Psi$, or
the last two rows are linearly dependent, when considered as
linear forms on $V_x\times V_y$ (for generic $z$).\par

If the first two rows are linearly dependent for every $(x,y,z)$,
putting $v:=F(x,y,z)\in V_{\xy}$, we see that there is some
coefficient vector
$(\alpha(x,y,v),\beta(x,y,v))\in\RR^2\setminus\{0\}$ such that
\begin{eqnarray}\label{8ad}
(\alpha(x,y,v)Av+\beta(x,y,v)Bv)\cdot\xi=0\quad \text{for all}\ 
\xi\in V_x\, ,
\end{eqnarray}

\noi for every  $\xy\in U ,v\in V_{\xy} .$

Choosing $\xi=v\in V_{\xy}$, in view of \eqref{8ab}  this yields
\[
\alpha(x,y,v)Q_A(v)=0,
\]
and by Lemma \ref{9'} (compare also \eqref{8tt}), since $r\ge 9$,
$Q_A(v)\not= 0$ for generic
$v$, so that $\alpha(x,y,v)=0$ for generic $v\in V_{\xy}$, hence
$\beta(x,y,v)\not= 0$ for generic $v\in V_{\xy}$. Thus
\[
(Bv)\cdot\xi =0\;\mbox{for generic}\  v\in V_{\xy}\;,\mbox{and
all}\,\xi\in V_x\, .
\]
But then this holds for all $v\in V_{\xy},\,\xi\in V_x$, so that
\begin{equation}\label{8ae}
B(V_{\xy})\subset\, \span_\RR \{Ax,Bx\}\quad \text{for all}\ \xy\in
U\, .
\end{equation}
We now distinguish two cases.

\medskip

\noi {\bf Case a.1:} $\rank \Psi =2$ for every $\xy\in U$ and
generic
$z$.
\medskip

 Then the first two lines, and also the last two lines of $\Psi$
are always linearly dependent, for generic $z$, and the preceding
discussion, in particular \eqref{8ae}, shows that
\[
B(V_{\xy})\subset\, \span_\RR\{Ax,Bx\}\, ,
\]
and analogously also
\[
B(V_{\xy})\subset\, \span_\RR\{Ay,By\}\, ,
\]
for all $\xy\in U.$ Thus
\[
B(V_{\xy})\subset\, \span_\RR\,\{Ax,Bx\}\cap\,
\span_\RR\{Ay,By\}=\{0\}\, ,
\]
hence $B|_{V_{\xy}}=0$. Since $\bigcup_{\xy\in U}
V_{\xy}$ spans $\RR^n$, this implies again $B=0.$

\noi There remains

\medskip

\noi {\bf Case a.2:} rank $\Psi=3$ for generic $(x,y,z)\in
U\times\RR^{n-4}$.

\medskip
 Shrinking $U$, if necessary, we can then assume
that, e.g., the last two rows of $\Psi$ are linearly independent
for generic $(x,y,z)$, and consequently the first two rows are
linearly dependent, for every $(x,y,z)\in U\times\RR^{n-4}$.
Freezing $x$, for generic $x$, and applying the same reasoning as
before to the mapping $F_x:(y,z)\mapsto F(x,y,z)$ instead
of $(x,y,z)\mapsto F(x,y,z)$ (by setting $\xi=0$), we see that rank
${\psi_3\choose\psi_4}=2$ implies that $D_{(y,z)}F_x(x,y,z)$ has
rank $n-2$, for generic $(y,z)$. This implies that, for every $x$
in some open set $U_1\subset \N_0$, the image $W_x$ of $F_x$
contains an analytic submanifold $\Om_x$ of dimension $\ge
n-2$.\par

But, if $U_2$ is an open subset of $\N_0$ such that $U_1\times
U_2\subset U$, then
\begin{eqnarray*}
\Om_x\subset \W_x  :=  \{F(x,y,z):y\in U_2,z\in\RR^{n-4}\}
                  =   \bigcup_{y\in U_2} V_{\xy}\subset V_x\, .
\end{eqnarray*}
Since $\dim\,V_x=n-2$, we see that $\Om_x$ is an open subset of
$V_x$.\par

And, since \eqref{8ae} still holds in the present case, we see that
$B(\Om_x)\subset\, \span_\RR\{Ax,Bx\}$, hence
\begin{equation}\label{8af}
B(V_x)\subset\, \span_\RR\{Ax,Bx\}\quad \text{for all}\  x\in U_1\,
.
\end{equation}
By Lemma \ref{8XX}, we can find non-empty open subsets
$\tilde U_1,\tilde U_2$ in $U_1$ such that $U':=\tilde U_1\times
\tilde U_2\subset\G$, and so that
$\bigcup_{\xy\in U'}V_{\xy}$ spans $\RR^n$. Then, by \eqref{8af},
also
\[
B(V_y)\subset\, \span_\RR\{Ay,By\}\quad \text{for all}\  y\in U_2\,
,
\]
so that
\begin{eqnarray*}
B(V_{\xy}) & = & B(V_x\cap V_y)\subset B(V_x)\cap B(V_y)\\
           & \subset & \span_\RR\{Ax,Bx\}\cap\,
\span_\RR\{Ay,By\}=\{0\},
\end{eqnarray*}
for all $\xy\in U'$. Again, we arrive at the contradiction
that
$B=0$.
\medskip

\noi {\bf (b)} Assume next that $m=3$.
\medskip

\noi Except for  an exchange of the r\^oles of $x$ and $y$, there
are then the following two possibilities: 
\smallskip

\noi {\bf (b.1)} $Ax\wedge Bx\wedge Ay\not= 0$, for some
$\xy=(x_0,y_0)\in U.$

Then, after shrinking $U$, if necessary, the same holds  for
all $\xy\in U,$ so that 
\[ N_{\xy}=\, \span_\RR\{Ax,Bx,Ay\}\quad \text{for all}\ \xy\in U\,
.
\]

\noi {\bf (b.2)} $Ax\wedge Bx\wedge By\not= 0$, say, again, for
all
$\xy\in U$ (after shrinking $U$). Then
\[
N_{\xy}=\,\span_\RR\{Ax,Bx,By\}\quad  \text{for all}\ \xy\in U.
\]

We begin with Case (b.1).
Since the arguments are quite similar to the ones in Case (a), we
shall content ourselves with a brief sketch, just indicating the
necessary modifications.

Choose again an orthonormal basis
\[
e_1\xy,\ldots ,e_{n-3}\xy
\]
of $V_{\xy}$, varying analytically in $\xy\in U$, and put, for
$(x,y,z)$ fixed,
\[
\Psi (\xi ,\eta)=\left(\begin{array}{ccc}
\psi_1(\xi,\eta)\\
\psi_2(\xi ,\eta)\\
\psi_3(\xi ,\eta)\end{array} \right) =-\left(\begin{array}{cc}
                                              \trans (A\cdot F(x,y,z)) & 0\\
                                              \trans (B\cdot F(x,y,z)) & 0\\
                                             0                    &  \trans (A\cdot F(x,y,z))\\
                                              \end{array} \right)
                                             {\xi \choose \eta}\,
                                             .
\]
Then, for $\xy\in U$,
\[
2\le\, \rank\Psi\le 3,\quad\mbox{for generic}\, z .
\]

If $\rank\Psi=3$ for some $(x,y,z)$, then the image  of $F$ contains
again an open subset of $\RR^n$, and we conclude again that
$B=0$.\par

So, assume that $\rank\Psi=2$, say, for all $\xy\in U$,
$z\not=0$. Then, the first two rows of $\Psi$ are linearly
dependent, so that
\eqref{8ae} holds. Moreover, similarly as in Case (a.2), for
fixed $x$, the mapping
$F_x:(y,z)\mapsto F(x,y,z)$ contains an analytic submanifold of
dimension $\ge (n-3)+1=n-2$ in its image, which is itself
contained in $V_x$. Since $\dim\,V_x=n-2$, we can conclude as in
Case (a.2).
\medskip

We are left with Case (b.2). Here,
\[
Ay\in\,\span_\RR\{Ax,Bx,By\}\quad \text{for all}\ \xy\in U\, .
\]
Let us assume that $U$ is the direct product $U=U_1\times U_2$ of
open sets $U_1,U_2\subset\N_0$. Then
\begin{equation}\label{8ag}
Ay\in N_x+\range B\quad \text{for all}\  y\in U_2\, ,
\end{equation}
for every $x\in U_1$.\par Now, since $Q_{B,x,y}=0$ and $\dim\,
V_{\xy}=n-3$, we see that $\rank\, B\le 2\cdot 3=6$, so that
\[
\dim\, (N_x+\range B)\le 8\, .
\]
On the other hand, since $U_2$ spans $\RR^n,$ and since $\rank A=r,$
\eqref{8ag} implies that
\[
\dim\, (N_x+\range B)\ge r .
\]
In combination, we find that $r\le 8$, contradicting our
assumptions.

 \qed

Combining  Proposition \ref{8pp} and  Lemma \ref{8SS} to
Lemma \ref{8YY}, we obtain

\begin{theorem}\label{MF2}
Let $A,B,C\in \Sym(n,\RR)$ satisfy our Standing Assumptions
\ref{ass}. Assume that $\maxrank\{A,B\}\ge 17,$   and that at least
one connected component of $\N$ spans
$\RR^n.$ Then $C$ is a linear combination of $A$ and $B$.
\end{theorem}

Applying the same type of technics, we can now obtain further
information also on the case where no  stratum of $\N$ spans.

\begin{lemma}\label{rank2}
Let $A,B\in \Sym(n,\RR)$ form  a non-dissipative pair, and
assume  that \hfill\newline 
$\maxrank\{A,B\}\ge 9.$ Then the following are equivalent:

\be 
\item[(i)]  There is a connected component of
$\N$ which does not span $\RR^n.$
\item[(ii)]No connected component of
$\N$   spans $\RR^n.$
\item[(iii)] $\minrank\{A,B\}=2.$
\ee
 
\end{lemma}

\noi {\bf Proof}. (iii) $\implies $ (ii). If $\minrank\{A,B\}=2,$
then we may assume without loss of generality that $\rank B=2.$ This
means that there are linearly independent vectors $\xi,\eta\in
\RR^n$ such that
$Q_B(x)=(\xi\cdot x)(\eta\cdot x).$ But then clearly every
component of $\N$ lies in one of the subspaces $\xi^\perp,$ or
$\eta\perp.$

\noi (iii) $\implies $ (ii) is trivial.

\noi (i) $\implies $ (iii).
Assume that there is a
connected component
$\N_0$ of $\N$ which lies in a subspace $\nu^\perp,$ where $\nu$ is
a unit vector. Without loss of generality, we may also assume that
$\rank A=
\maxrank\{A,B\}\ge 9. $ Let $x_0\in\N_0.$ Arguing as in the
proof of Lemma
\ref{8JJ}, we see that 
\begin{equation}\label{r21}
\nu\cdot x=\al(x)Q_A(x)+\beta(x) Q_B(x), 
\end{equation}
for all $x$ in a sufficiently small neighborhood of $x_0.$ Here,
$\al$ and $\beta$ are analytic functions near $x_0.$ Applying
the second derivative, and restricting the forms to $V_{x_0},$ we
obtain
\begin{equation}\label{r22}
0= \alpha (x_0)Q_{A,x_0}+\beta(x_0) Q_{B,x_0}\quad \text{for all}\
x_0\in\N_0.
\end{equation}
Notice that $\beta (x_0)\ne 0,$ since $Q_{A,x_0}\ne 0,$ in view of
Lemma \ref{9'}.

Applying next the same kind of reasoning as in the proof of Lemma
\ref{8SS}, only with $V_{(x,y)}$ replaced by $V_x$ and $ Q_{A,x,y}$
by $Q_{A,x},$ etc., and $m=2,$ we see that \eqref{r22} implies that
there is a non-trivial linear combination $D=\al_0 A+\beta_0B$ such
that 
\begin{equation}\label{r23}
Q_{D,x}=0\quad \text{for all}\ x\in\N_0.
\end{equation}

Notice also that 
\begin{equation}\label{r24}
\span_\RR \N_0=\nu^\perp.
\end{equation}
Indeed, otherwise $\N_0$ would be an open subset of  a linear
subspace $W$ of dimension $n-2,$ and $Q_A|_W=0.$ But this would
imply $\rank A\le 4,$ contradicting our assumption on $A.$

Arguing similarly as in the proof of LemmaÊ\,\ref{8XX}, \eqref{r24}
implies that 
\begin{equation}\label{r25}
\bigcup_{x\in U}V_x \quad \text{spans}\ \nu^\perp,
\end{equation}
for every non-empty open subset $U$ of $\N_0.$ Indeed, otherwise
$\bigcup_{x\in U}V_x $ would be contained in a  linear
subspace $W$ of dimension $n-2,$ so that, by comparing dimensions,
$V_x=W$ for every $x\in U.$ But this would imply $\N_0\subset W,$
contradicting \eqref{r24}.

Finally, we can apply a similar reasoning as in the proof of Lemma
\ref{8YY} in order to conclude that 
\begin{equation}\label{r26}
Q_D|_{\nu^\perp}=0.
\end{equation}
By Lemma \ref{9'}, this implies $\rank D\le 2,$ hence $\rank D=2,$
and thus $\minrank\{A,B\}=2.$

To prove \eqref{r26}, we may assume that there is an orthonormal
basis
$e_1(x),\ldots ,e_{n-2}(x)$ of $V_{x}$, varying analytically in
$x\in U$. We then put 
$$f_1(x):=Ax, \ f_2(x):=Bx,$$
 so that 
\[
e_1(x),\ldots ,e_{n-2}(x),f_1(x),f_2(x)
\]
is a basis of $\RR^n$. Consider the mapping
$F:U\times\RR^{n-2}\rightarrow \nu^\perp,$
\[
F(x,z):=\sum^{n-2}_{j=1}z_je_j(x)\in V_x\subset \nu^\perp,\quad x\in
U,z\in\RR^{n-2}\, , \]
where we consider $U$ as an analytic submanifold of
$\RR^{n}$ of dimension $n-2.$
\medskip

We shall prove that there is some $(x_0,z_0)$ such that 
\begin{equation}\label{r27}
\rank\, DF (x_0,z_0)=n-1\, .
\end{equation}
This implies that $F$ is a submersion near $(x_0,z_0)$, so
that $\bigcup_{x\in U}V_x$  contains a non-empty
open subset $\Om$ of $\nu^\perp.$ Since, by \eqref{r23}, $Q_B$
vanishes on $\Om$, we thus obtain \eqref{r26}. 

In order to prove \eqref{r27}, we consider here, for $x,z$ fixed,
the linear mappings 
\begin{eqnarray*}
\psi_i &:=& V_x \to \RR\, ,\\
\psi_i(\xi)&:=& f_i(x)\cdot (D_{x}F(x,z)(\xi)),\quad
i=1,2\, , 
\end{eqnarray*}
and 
\[
\Psi:=\left(\begin{array}{c} \psi_1\\ 
\psi_2\end{array}\right) :V_x\rightarrow\RR^2. 
\]
Similarly as in the proof of Lemma \ref{8YY}, it then suffices to
show that linear mapping $\Psi$ 
has rank 1, generically. But, 
\begin{eqnarray*}
\Psi (\xi)=
-\left(\begin{array}{c} \trans (A\cdot F(x,z)) \\
 \trans (B\cdot F(x,z))  \end{array} \right)
\xi \, ,
    \end{eqnarray*}
and, for $x\in U,$ there is some
$z\in\RR^{n-2}$ such that $\psi_1|_{V_{x}}\not= 0$. Indeed,
otherwise we would have 
$(Av)\cdot\xi =0\quad \text{for all}\  v\in V_{x}$, $\xi\in V_x$,
and in particular $Q_{A,x}=0$. But then rank $A\le 2\cdot 2=4$,
 a contradiction. This shows that $\psi_1\not= 0,$ hence $\rank
\Psi= 1.$

\qed

\medskip

Our  main result Theorem \ref{8ZZ}  concerning the form problem is now an 
immediate consequence of  Theorem \ref{8Z}, 
Theorem \ref{MF2} and Lemma \ref{rank2}.

\begin{remark}\label{genform}
The statement of Theorem \ref{8ZZ} remains true, if we replace the
canonical symplectic form on $\RR^{2d}$ by an arbitrary constant symplectic
form $\om$,  and define the Poisson bracket accordingly by
\[
\{f,g\}_\om :=\om (H^\om_f,H^\om_g)\, ,
\]
where $H^\om_f$ denotes  the Hamiltonian vector field associated to $f$ (see
Section \ref{introduction}), 
 and define
$C$ by
$$Q_C:=\{Q_A,Q_B\}_\om.$$
\end{remark}

To see this, notice that there is a 
linear change of coordinates which transforms
$\om$ into the canonical symplectic form on $\RR^{2d},$ and that 
the space of quadratic forms and the cone of positive-semidefinite
forms remain invariant under such a change of coordinates.

\setcounter{equation}{0}
\section{Applications to non-solvability of doubly characteristic differential operators}

We are now in a position to prove Theorem \ref{9E}.
Assume that $L$ satisfies the assumptions of this theorem. Recall that $q_j$ denotes
the symbol of $X_j,$ so that 
$f(x,\xi):={\dst\sum^m_{j,k=1}}a_{jk}(x)q_j(x,\xi)q_k (x,\xi)$ and
$g(x,\xi):={\dst\sum^m_{j,k=1}}b_{jk}(x)q_j(x,\xi)q_k (x,\xi)$
denote the real and imaginary parts of the principal symbol of
$L$.

According to H\"ormander's Theorem \ref{Ho},  it suffices to prove
that there is some $\eta\in\RR^n$ such that
\begin{equation}\label{9f}
f(x_0,\eta)=g(x_0,\eta)=0\ \mathrm{and}\ \{f,g\}(x_0,\eta)\not=0\,
.
\end{equation}
One easily computes that
\[
\{f,g\}(x,\xi)=4\;\sum_{j,k}(A(x)J_{(x,\xi'')}B(x))_{jk}\,
q_j(x,\xi)q_k(x,\xi)+h(x,\xi)\, ,
\]
where $h$ is of the form
\begin{equation}\label{9g}
h=\sum_{j,k,l}d_{jkl}\, q_j q_h q_l\, ,
\end{equation}
with smooth functions $d_{jhl}=d_{jhl}(x)$ depending only on
$x$. Since, as one easily checks,
\[
\{Q_{A(x)},Q_{B(x)}\}_{(x,\xi'')}=4 \trans v A(x)J_{(x,\xi'')}
B(x)v,\quad v\in\RR^m\, ,
\]
we thus get

\begin{equation}\label{9h}
\{f,g\}(x,\xi)=Q_{C(x,\xi'')}\left(q(x,\xi)\right)\,+\, h(x,\xi)\, ,
\end{equation}
if we put $q(x,\xi):=(q_1(x,\xi),\ldots ,q_m(x,\xi))$.
\smallskip

Putting $\xi'_0:=-E(x)\,\xi''_0,$ we have $q(x_0,\xi_0)=0$, i.e.,
$(x_0,\xi_ 0)\in\Sigma_2$, so that
\[
q(x_0,\xi_0+(v',0))=v'\quad \text{for every}\  v'\in\RR^m\, .
\]
But, by our assumptions, Theorem \ref{8ZZ} (more precisely, Remark
\ref{genform} ) implies that there is some $v'\in\RR^m$ such that
\[
Q_{A(x_0)}(v')=Q_{B(x_0)}(v')=0\ \mathrm{and}\ 
Q_{C(x_0,\xi''_0)}(v')\not= 0\, .
\]
For $0<t\le 1$, let $\xi(t):=\xi_0+(tv',0)$. Then
$q(x_0,\xi(t))=tv'$, so that
\[
\begin{array}{lcl}
 f(x_0,\xi(t)) & = & t^2Q_{A(x_0)}(v')=0\, , \\
 g(x_0,\xi(t)) & = & t^2 Q_{B(x_0)}(v')=0 \, ,\\
\end{array}
\]
and, by \eqref{9g}, \eqref{9h},
\[
\{f,g\}(x_0,\xi(t))=t^2 Q_{C(x_0,\xi(t)'')} (v')+O(t^3)\, .
\]
But, as $t\rightarrow 0$, $C(x_0,\xi(t)'')\rightarrow
C(x_0,\xi''_0)$, so that 
$$\dst{\lim_{t\rightarrow 0}}\,
t^{-2}\{f,g\}(x_0,\xi(t)'')=Q_{C(x_0,\xi''_0)}(v')\not= 0.$$
For sufficiently small $t_0>0$, we thus have
\[
\{f,g\}(x_0,\xi(t_0))\not=0\, ,
\]
so that \eqref{9f} is satisfied for $\eta :=\xi_0+(t_0v',0)$.

\qed
\medskip

As an application of Theorem \ref{9E},  consider a connected, simply connected
two-step nilpotent Lie group $G$. Up to an automorphism, such a group can
always be realized as $G=\RR^m\times\RR^{\ell}$ (as a manifold),
with a group law of the form

\[
(z,u)\cdot (z',u')=(z+z',u+u'+\frac 1 2( \trans z J^{(i)}z')_i)\, ,
\]
for $z,z'\in\RR^m, u,u'\in\RR^{\ell}$, where $J^{(1)},\ldots
,J^{(\ell)}$ are skew-symmetric $m\times m$-matrices. A basis of
the Lie algebra ${\frak g}$ of left-invariant vector fields is
then given by

\begin{eqnarray*}
  X_j & := & \frac{\de}{\de z_j}+\frac 1 2 \sum^{\ell}_{i=1}(
\trans z\cdot J^{(i)})_j \,\frac{\de}{\de u_i},\qquad  j=1,\ldots
,m\, , \\
  U_k & := & \frac{\de}{\de u_k}, \qquad  k=1,\ldots ,\ell\,
. 
\end{eqnarray*}

Consider an operator $L$ of the form

\begin{equation}\label{9i}
L=\sum^m_{j,k=1}\al_{jk} X_j X_k +\;\mbox{lower order terms}\, ,
\end{equation}
where the coefficient matrix ${\A}=(\al_{jh})_{j,h}\in\Sym (m,\CC)$
is symmetric. We put  $A:=\Re{\A},$ 
$B:=\Im{\A}$. The matrix $\{E_{jl}(x)\}_{jl}$ in \eqref{9d} is
then given by
\[
E_{ji}(z,u)=\frac 1 2\sum^m_{k=1} J^{(i)}_{kj} z_k\, ,
\]
and, if we split $\xi=(\zeta,\mu)\in\RR^m\times\RR^{\ell}$ according
to the coordinates $x=(z,u)\in G$, we have

\[
J_{((z,u),\mu)}=J^{\mu}\, ,
\]
if we put $J^{\mu}:={\dst\sum^{\ell}_{i=1}}\mu_i J^{(i)}$.\par

Finally, we define 
\[
C_{\mu_0}:=\frac 1 2 (AJ^{\mu_0}B-BJ^{\mu_0}A)\, .
\]
From Theorem \ref{9E}, we then immediately obtain

\begin{cor}\label{9j}
Assume that $J^{\mu_0}$ is non-degenerate for some
$\mu_0\in\RR^{\ell}$, that $A,B$ forms a non-dissipative pair, and 
that $A,B$ and $C_{\mu_0}$ are linearly independent. Moreover, 
suppose that either 
\be
\item[(i)]   $\minrank \{A, B\}\ge 3$ and 
$\maxrank\{A, B\}\ge 17, $ or
\item[(ii)] $\minrank \{A, B\}= 2,$
$\maxrank\{A, B\}\ge 9, $ and that the joint kernel 
 $\ker A\cap \ker B $ of $Q_A$ and $ Q_B$ is
either trivial, or a symplectic subspace with respect to the
symplectic form $\om_{\mu_0}$ on $\RR^m$ associated to $\trans
(J^{\mu_0})^{-1}.$
\ee 

\noi Then the operator $L$ in \eqref{9i} on $G$ is
nowhere locally solvable.
\end{cor}

In the special case of the Heisenberg group $\HH_d$, we have
$m=2d,\ \ell=1$ and $J^{(1)}:=J=\left(%
\begin{array}{cc}
  0 & I_d \\
  -I_d & 0 \\
\end{array}%
\right)$.  Then  $J_{\mu}=\mu J$, $\mu\in\RR,$ so that we may choose $\mu=1$. Putting
\[
C:=\frac 1 2 (AJB-BJA)\, ,
\]
we obtain

\begin{cor}\label{9k}
If $G=\HH_d$ is  the Heisenberg group, assume 
 that $A,B$ forms a non-dissipative pair, and 
that $A,B$ and $C$ are linearly independent. Moreover, 
suppose that either 
\be
\item[(i)]   $\minrank \{A, B\}\ge 3$ and 
$\maxrank\{A, B\}\ge 17, $ or
\item[(ii)] $\minrank \{A, B\}= 2,$
$\maxrank\{A, B\}\ge 9, $ and that the joint kernel 
 $\ker A\cap \ker B $ of $Q_A$ and $ Q_B$ is
either trivial, or a symplectic subspace with respect to the
canonical symplectic form on $\RR^{2d}$ (associated to $-J$).
\ee 

\noi Then the operator $L$ in \eqref{9i} on $\HH_d$ is
nowhere locally solvable.
\end{cor}

This shows that local solvability of $L$ on  Heisenberg groups can
essentially only arize  if the operator 
$e^{i\theta}L$ is dissipative, for some $\theta\in\RR$.  This
statement is true in the strict sense, if, e.g., the matrix $A$ is
non-degenerate, and $d\ge 9.$ \par

As we had already mentioned, the examples in [KM] and [MP] show
that the analogous statement is wrong on Heisenberg groups of low
dimension 5 and 7. \\

\begin{remark}\label{3step}
Theorem \ref{9E} applies also to higher
step situations.

\end{remark}

 For instance, assume that ${\g}$ is a nilpotent Lie
algebra of step $r\ge 3$, and let
$\g=\g_1\subset\g_2\subset\cdots\subset\g_{r+1}=\{0\}$ denote the
descending central series, i.e., $\g_{j+1}:=[\g ,\g_j]$. Let
$G=\exp\g$ be the associated nilpotent Lie group, and choose
elements $X_1,\ldots ,X_m$ of $\g$ which form a basis modulo
$\g_2$. Consider the $X_j$ as left-invariant vector fields on $G$
as usual, and let $L$ on $G$ be given by \eqref{9i}. Put
$G_j:=\exp\g_j,$ and let $H:=G_3\setminus G=\{G_3g:g\in G\}$ denote
the quotient group of $G$ by $G_3.$ Then
$H$ is a 2-step nilpotent Lie group, and if we assume that also
the lower order terms in \eqref{9i} are left-invariant, then $L$
factors through $H$ as a left-invariant differential operator

\[
\tilde{L}=\sum^m_{j,k=1}\al_{jk}\tilde{X}_j\tilde{X}_k+\;\mbox{lower
order terms}\, .
\]
Here, $\tilde{X}_j$ is the left-invariant vector field on $H$
corresponding to $X_j$.\par

Choose further elements $X_{m+1},\ldots ,X_N$ such that
$X_1,\ldots ,X_N$ forms a basis modulo $\g_3$, and then elements
$Y_1,\ldots ,Y_k$ so that $X_1,\ldots ,X_N,Y_1,\ldots ,Y_k$ forms
a basis of $\g$. We may then choose coordinates
$(x,y)\in\RR^N\times\RR^k$ of $G$, by putting
\[
g(x,y):=\exp\left( \sum^k_{j=1}y_jY_j\right)\exp \left(
\sum^N_{\ell=1}x_{\ell}X_{\ell}\right)\, .
\]

Let $(\xi ,\eta)$ denote the dual variables. By  $2\pi \sigma (A)$
we denote the principal symbol of a
differential operator $A.$ Then one easily shows that 
\[
\sigma (X_j)((x,0),(\xi ,0))=\sigma (\tilde{X}_j)(x,\xi)\, ,
\]
hence
\begin{equation}\label{9l}
\sigma (L)((x,0),(\xi ,0))=\sigma (\tilde{L})(x,\xi)\quad
 \text{for all}\ \xi\in\RR^N
\end{equation}
(here, we have chosen $x$ as natural coordinates for $H$).
Moreover, since
\[
[X_j,X_k]^{\sim}=[\tilde{X}_j,\tilde{X}_k]\, ,
\]
we have
\begin{eqnarray*}
\{\sigma(X_j),\sigma(X_k)\}((x,0),(\xi
,0))=i\sigma([X_j,X_k])((x,0), (\xi ,0))\\
=i\sigma([X_j,X_k]^{\sim})(x,\xi)=\{\sigma(\tilde{X}_j),\sigma(\tilde{X}_k)\}(x,\xi)\,
.
\end{eqnarray*}

If $\tilde{J}_{(x,\xi)}$ denotes the skew symmetric matrix on $H$
given in Theorem \ref{9E}, and $J_{((x,0),(\xi ,0))}$ the one on
$G$, we thus have

\[
J_{((x,0),(\xi ,0))}=\tilde{J}_{(x,\xi)}\, .
\]

Thus, if we assume that $\tilde{L}$ satisfies the hypotheses of
Corollary \ref{9j}, then with $\tilde{L}$, also $L$ satisfies the
assumptions of Theorem \ref{9E}, so that $L$ is nowhere locally
solvable.

\bibliographystyle{plain}
%\bibliography{ReA>0}
%\end{document}

\vskip1cm

\noindent\emph{Mathematisches Seminar,  C.A.-Universit\"at Kiel,
Ludewig-Meyn-Str.4, D-24098 Kiel, Germany\\
 e-mail: mueller@math.uni-kiel.de}
\end{document}